\documentclass[11pt]{article}

\usepackage[utf8]{inputenc}
\usepackage{amsmath,mathtools}
\usepackage{amsfonts,amssymb,latexsym,multicol,color}
\usepackage[francais]{babel}
\usepackage[T1]{fontenc}
\usepackage{graphicx}

\newcounter{fig}
\input{cyracc.def}

\font\tencyrit=wncyi10 scaled 1100

\font\tencyr=wncyr10 scaled 1000
\def\cyr{\tencyr\cyracc}

\def\cyrit{\tencyrit\cyracc}

\newtheorem{ques}{Question}
\newtheorem{theo}{Th\'eor\`eme}

\newtheorem{prop}{Proposition}

\newcommand{\cad}{\text{c'est-\`a-dire }}
\newcommand{\expli}[1]{\quad\text{\footnotesize (#1)}}

\makeatletter
\ifcase\@ptsize

\or

\or

\fi
\makeatother

\newcommand{\implique}{\Rightarrow}

\newcommand{\ioe}{\leqslant}
\newcommand{\soe}{\geqslant}

\newcommand{\vers}{\rightarrow}
\newcommand{\dist}{{\rm dist}}

\newcommand{\demi}{{\frac{1}{2}}}

\newcommand{\mes}{{\rm mes}}

\newcommand{\sgn}{{\rm sgn}\,}

\newcommand{\eps}{\varepsilon}
\newcommand{\fhi}{\varphi}
\newcommand{\thet}{\vartheta}

\newcommand{\abs}[1]{\left\lvert #1 \right\rvert}
\newcommand{\floor}[1]{\left\lfloor #1 \right\rfloor}
\newcommand{\norm}[1]{\left\lVert #1 \right\rVert}
\newcommand{\set}[1]{\left\{ #1 \right\}}

\newcommand{\Ccal}{{\mathcal C}}

\newcommand{\Ecal}{{\mathcal E}}
\newcommand{\Fcal}{{\mathcal F}}

\newcommand{\Kcal}{{\mathcal K}}
\newcommand{\Lcal}{{\mathcal L}}

\newcommand{\Pcal}{{\mathcal P}}

\newcommand{\Tcal}{{\mathcal T}}

\newcommand{\Nat}{{\mathbb N}}
\newcommand{\Int}{{\mathbb Z}}
\newcommand{\Rat}{{\mathbb Q}}
\newcommand{\Real}{{\mathbb R}}

\newcommand{\fin}{\hfill$\Box$}
\newcommand{\dem}{\noindent {\bf D\'emonstration\ }}
\newcommand{\fine}{\tag*{\mbox{$\Box$}}}

\providecommand{\bysame}{\leavevmode ---\ }
\providecommand{\og}{``}
\providecommand{\fg}{''}
\providecommand{\smfandname}{et}
\providecommand{\smfedsname}{\'eds.}
\providecommand{\smfedname}{\'ed.}

\newcommand{\virg}{\raisebox{.7mm}{,}}

\title{Autour d'un problème extrémal étudié par Edmund Landau}

\author{Michel Balazard}

\date{}

\begin{document}
\maketitle

\begin{center}
  {\sc Abstract}
\end{center}
\begin{quote}
{\footnotesize This expository article is an introduction to Landau's problem of bounding the derivative, knowing bounds for the function and its second derivative, and some of its variants and generalizations. Connexions with convex and functional analysis, numerical integration, and approximation theory are emphasized. Among others, we describe the set of extreme points of the relevant set of functions.}
\end{quote}

\begin{center}
  {\sc Keywords}
\end{center}
\begin{quote}
{\footnotesize Extremal problems, Landau inequality, Kolmogorov inequality, extreme points, Peano kernel, Euler splines. \\MSC classification : 41A17}
\end{quote}


\section{Introduction}

Soit $I$ un intervalle de la droite réelle $\Real$, et $n$ un nombre entier supérieur ou égal à $2$. Nous allons considérer des fonctions $f : I \vers \Real$ vérifiant les propriétés suivantes,

\smallskip

$(i)$ $f$ est $(n-1)$ fois dérivable sur $I$ ;

$(ii)$ la dérivée $f^{(n-1)}$ est absolument continue sur tout segment inclus dans $I$.

\smallskip

Dans $(i)$ la dérivabilité est à comprendre à droite (resp. à gauche) en la borne inférieure (resp. supérieure) de $I$, si celle-ci appartient à $I$. L'hypothèse $(ii)$ d'absolue continuité de $f^{(n-1)}$ entraîne l'existence de $f^{(n)}(x)$  pour presque tout $x \in I$. L'ensemble des fonctions $f$ vérifiant~$(i)$ et $(ii)$ est un espace vectoriel réel que nous noterons $V_n(I)$. Observons que $V_n(I) \subset V_m(I)$ si~$2 \ioe m\ioe n$.

\smallskip

Si $a$ et $b$ sont deux nombres réels positifs, soit $\Lcal_n(a,b\, ; \, I)$ l'ensemble des fonctions $f \in V_n(I)$ vérifiant, outre $(i)$ et $(ii)$, les inégalités suivantes,

\smallskip

$(iii)$ on a $\lvert f(t) \rvert \ioe a $ pour tout $t \in I$ ;

$(iv)$ on a $\lvert f^{(n)}(t) \rvert \ioe b$ pour presque tout $t \in I$.

 \smallskip

La conjonction de $(ii)$ et $(iv)$ est équivalente à l'hypothèse que $f^{(n-1)}$ est lipschitzienne de rapport $b$ :
\[
\lvert f^{(n-1)}(t_1) - f^{(n-1)}(t_2) \rvert \ioe b \lvert t_1-t_2\rvert \quad (t_1,t_2 \in I).
\]

\`A la suite de Schoenberg (cf. \cite{zbMATH03411359}, \S 5, p. 131), Chui et Smith ont signalé dans \cite{zbMATH03501820} l'interprétation cinématique de cette définition pour le cas $n=2$ : les fonctions~$f$ de $\Lcal_2(a,b\, ; \, I)$ sont les lois des mouvements de particules contraintes à rester dans l'intervalle~$[-a,a]$ de la droite, pendant l'intervalle de temps~$I$, et dont l'accélération est majorée par~$b$ en valeur absolue. J'emploierai parfois ce langage.

Dans son article \cite{zbMATH02621634} de 1914, Landau a posé, et essentiellement résolu, le problème de déterminer quelle vitesse instantanée maximale une telle particule peut atteindre. La solution complète est explicitée dans l'article de Chui et Smith mentionné ci-dessus. 

Ce travail de Landau, motivé par des considérations de théorie taubérienne, alors en plein développement (cf. \cite{zbMATH02625426}, II, p. 416-428), a ensuite suscité une intense activité dans l'étude des majorations des dérivées d'ordres intermédiaires $f^{(k)}$, $0<k<n$, lorsque $f \in \Lcal_n(a,b\, ; \, I)$. Appelons \emph{problème de Landau} la question de déterminer les bornes supérieures 
\begin{equation}\label{200425b}
\sup_{f \in \Lcal_n(a,b\, ; \, I)} \; \sup_{t \in I}\, \lvert f^{(k)}(t) \rvert \quad (0<k<n),
\end{equation}
et, le cas échéant, de préciser les fonctions~$f$ extrémales.

Le résultat central de la théorie est la solution complète de ce problème en~1939 par Kolmogorov (cf. \cite{zbMATH00050290}, pp. 277-290), dans le cas où $I=\Real$.

Le présent article est une introduction à cette théorie. J'ai choisi de centrer l'exposé sur le cas~\mbox{$n=2$}, qui correspond au problème initial de Landau et permet de présenter les idées, également pertinentes pour le cas général, sous leur forme la plus accessible, notamment lorsqu'il s'agira de variantes {\og ponctuelles\fg} et {\og en moyenne\fg} du problème de Landau. Cela étant, je décrirai pour $n$ quelconque certaines propriétés géométriques et topologiques de l'ensemble~$\Lcal_n(a,b\, ; \, I)$, et présenterai, le plus souvent sans démonstrations, les principaux résultats connus dans ce cas général.

Par ailleurs, je n'aborde pas systématiquement la question, parfois délicate, de décrire toutes les fonctions extrémales, solutions des divers problèmes extrémaux abordés.

\smallskip

Voici le plan de cet article. Au \S \ref{190322a}, l'invariance du problème de Landau par translation et dilatation est utilisée pour réduire son étude à trois cas particuliers. Nous portons notamment notre attention sur l'ensemble
\[
\Lcal_n(T)=\Lcal_n(1,1;\, [0,T]),
\]
où $T>0$. Le \S \ref{190322c} concerne la propriété géométrique fondamentale de l'ensemble $\Lcal_n(T)$, sa convexité, et contient une description de ses points extrêmes. Le \S\ref{190322b} présente succinctement la notion de \emph{noyau de Peano} ; dans le contexte du problème de Landau, elle permet de représenter les valeurs des dérivées intermédiaires $f^{(k)}$ comme fonctions linéaires explicites de $f$ et $f^{(n)}$. On en déduit une première majoration de $\abs{f^{(k)}}$, si $f \in \Lcal_n(T)$. Au \S \ref{190322d}, la compacité de $\Lcal_n(T)$ pour la topologie de la convergence uniforme des fonctions et de leurs dérivées d'ordre $<n$ permet d'appliquer le théorème de Krein et Milman. Le \S \ref{190322e} rassemble quelques généralités sur certains problèmes extrémaux dans l'ensemble $\Lcal_n(T)$. Le \S\ref{200416a} est dévolu à une étude approfondie du cas~$n=2$. La solution, très simple, du problème initial sur $\Lcal_2(T)$ est présentée au~\S\ref{190216a} ;  le cas de $I=[0,\infty[$ en découle au \S\ref{200423a}. Le \S\ref{190406a} contient la solution d'une variante ponctuelle du problème, cas particulier de l'étude exhaustive effectuée par Landau en~1925, dans son article~\cite{zbMATH02591405}. La notion de \emph{fonction de comparaison}, due à Kolmogorov, est introduite au \S\ref{200117b} ; elle fournit rapidement, au \S\ref{200416b}, la solution du problème de Landau pour~$I=\Real$. Après avoir établi trois propriétés de prolongement au \S\ref{190624b}, nous donnons, au~\S\ref{190406b}, les résultats de base sur le problème de Landau en moyenne, en utilisant les travaux de Bojanov et Naidenov (1999-2003, cf.~\cite{zbMATH01369964}, \cite{zbMATH01891483} et \cite{zbMATH01970534}). Enfin, le \S\ref{190528b} est un formulaire des résultats obtenus, en revenant au cas général de l'ensemble~$\Lcal_2(a,b\, ; \, I)$. Le \S\ref{200311a} est un survol des principaux résultats connus dans le cas général. Le \S\ref{200416c} présente les \emph{splines d'Euler}, fonctions extrémales pour le cas~$I=\Real$, comme l'a montré Kolmogorov ; son théorème est énoncé au \S\ref{200416d}. Le cas de l'intervalle~$I=[0,T]$ est l'objet du~\S\ref{200202a} ; celui de $I=[0,\infty[$ est présenté au \S\ref{200416e}. Enfin, le \S\ref{200416f} mentionne quelques résultats dans le cadre des espaces de Lebesgue $L^p$.

\smallskip

L'article contient \ref{190627e} propositions avec leurs démonstrations. Les plus importantes d'entre elles sont les propositions \ref{200118a}, \ref{200402a}, \ref{190605d}, \ref{190627g}, \ref{190605b}, \ref{190702a}, \ref{190704a}, \ref{190626b} et \ref{190615b}. Par ailleurs, deux parmi les nombreuses questions encore ouvertes de cette théorie sont mentionnées explicitement. 

\smallskip

Nous noterons 

$\bullet$ $\norm{f}_{\infty}$ la borne supérieure essentielle, éventuellement infinie, d'une fonction $f$ à valeurs réelles, supposée définie presque partout sur $I$, et mesurable ;

$\bullet$ $\floor{x}$ la partie entière, et $\set{x}$ la partie fractionnaire du nombre réel $x$, de sorte que
\[
x =\floor{x} + \set{x} \quad ; \quad \floor{x}  \in \Int \quad ; \quad 0 \ioe \set{x} <1.
\]

\section{Réduction de l'étude à trois cas}\label{190322a}

La restriction d'une fonction $f \in \Lcal_n(a,b\, ; \, I)$ à un sous-intervalle quelconque $J \subset I$ appartient à $\Lcal_n(a,b\, ;J)$.

Si $t_0$ est un nombre réel, point adhérent à $I$ et n'appartenant pas à $I$, la condition $\lvert f^{(n)} \rvert \ioe b$ presque partout sur $I$ entraîne que $f^{(n-1)}(t)$ admet une limite quand $t$ tend vers $t_0$. On en déduit que $f^{(j)} (t)$ admet également une limite $\ell_j$ quand $t$ tend vers $t_0$, pour tout $j$, $0 \ioe j <n$, et que la fonction $g$, prolongeant $f$ à $J=I \cup\{t_0\}$, telle que $g(t_0)=\ell_0$, appartient à $\Lcal_n(a,b\, ;J)$. Pour cette raison, nous supposerons~$I$ fermé dans toute la suite.

Ensuite, si $\lambda$ et $\mu$ sont des nombres réels non nuls, et $t_0 \in \Real$, la transformation $f \mapsto f_1$, où
\[
f_1(t)=\mu f\big(\lambda (t-t_0)\big),
\]
qui correspond à un changement d'unités de distance et de temps, et à un changement d'origine du temps, est une bijection entre $\Lcal_n(a,b\, ; \, I)$ et $\Lcal_n(a_1,b_1;I_1)$, où
\[
a_1=\abs{\mu}a \quad ; \quad b_1= \abs{\mu\lambda^n}b \quad ; \quad I_1= t_0+I/\lambda.
\]
En utilisant une telle transformation, on ramène l'étude de l'ensemble $\Lcal_n(a,b\, ; \, I)$ à l'un des trois cas suivants :
\[
\Lcal_n(1,1; \Real) \quad ; \quad \Lcal_n(1,1;[0,\infty[) \quad ; \quad \Lcal_n(1,1;[0,T]), \text{ où } T>0.
\]

C'est ce dernier cas, en particulier pour $n=2$, qui va nous occuper principalement dans cet article ; nous verrons notamment ses relations avec les deux premiers. Soit donc $T$ un nombre réel positif ; posons  
\[
I=[0,T] \quad ; \quad \Lcal_n(T)=\Lcal_n(1,1;I). 
\]
Les éléments de~$\Lcal_2(T)$ sont les lois des mouvements, de durée finie égale à $T$, d'une particule contrainte à rester dans le segment $[-1,1]$, et dont l'accélération est au plus~$1$ en valeur absolue. 

\smallskip

Notons toutefois que certains raisonnements sont présentés plus simplement en revenant au cas général de l'étude de $\Lcal_n(a,b\, ; \, I)$ (cf., par exemple, la proposition \ref{200402a}, au \S\ref{190322b}).

\section{Convexité}\label{190322c}

Le nombre positif $T$ étant fixé dans tout ce paragraphe, nous simplifions l'écriture en posant~$\Lcal_n=\Lcal_n(T)$ (et toujours $I=[0,T]$).

\smallskip

L'ensemble $\Lcal_n$ est une partie convexe de l'espace vectoriel $V_n(I)$ : 
\[
f,g \in \Lcal_n \text{ et } \vartheta \in [0,1] \implique \vartheta f+(1-\vartheta)g \in \Lcal_n.
\]

Rappelons que les \emph{points extrêmes} d'une partie $C$ d'un espace vectoriel réel sont les éléments~$x$ de $C$ tels que 
\[
\forall \, y,z \in C, \; \forall \, \vartheta \in\, ]0,1[\, , \; x=\vartheta y+(1-\vartheta)z  \implique y=z.
\]

Si $C$ est convexe, cette condition équivaut à
\[
\forall \, y,z \in C, \;  2x=y+z  \implique y=z.
\]

Nous déterminons maintenant les points extrêmes de $\Lcal_n$. La caractérisation présentée dans la proposition \ref{200118a} ci-dessous repose sur une notion de multiplicité adaptée à l'ensemble $\Lcal_n$. Si~$f \in \Lcal_n$, notons $F$ l'ensemble des $t\in I$ tels que $\abs{f(t)}=1$. Si~$t_0\in F$, la multiplicité de $t_0$ est, par définition, l'unique nombre entier $m$ tel que

$\bullet$ $1 \ioe m \ioe n$ ;

$\bullet$ $f^{(j)}(t_0)=0$ pour tout $j$ tel que $1\ioe j \ioe m-1$ ;

$\bullet$ $m=n$, ou $m<n$ et $f^{(m)}(t_0)\neq 0$.

Observons tout de suite que, si $t_0\in F\,\cap \, ]0,T[$ et si $m<n$, alors $m$ est pair, puisque $t_0$ est un point d'extremum de $f$.

\begin{prop}\label{200118a}

La fonction $f \in \Lcal_n$ est un point extrême de $\Lcal_n$ si, et seulement si les deux conditions suivantes sont vérifiées :

\begin{quote}
$(i)$ la somme des multiplicités des éléments du fermé $F=\set{t \in I, \, \abs{f(t)}=1}$ est supérieure ou égale à $n$ ;

$(ii)$ on a $\abs{f^{(n)} (t)}=1$ presque partout sur~$I\setminus~F$.

\end{quote}

\end{prop}
\dem

L'ensemble $F$ est bien un fermé pour toute fonction $f \in \Lcal_n$, puisqu'une telle fonction est continue. 

Nous allons adapter l'argumentation de Roy dans \cite{zbMATH03256144}, p. 1158 ; il y étudiait l'ensemble des fonctions lipschitziennes de rapport $1$ sur $[0,1]$. 

\smallskip

Pour commencer, montrons que les points extrêmes de $\Lcal_n$ vérifient $(i)$ et $(ii)$. 

Soit $f$ un tel point extrême. Supposons que $f$ ne vérifie pas $(i)$. L'ensemble $F$ est donc fini. Notons $t_i$, pour $1\ioe i \ioe k$, ses éléments deux à deux distincts, et $d_i$ la multiplicité de $t_i$ au sens ci-dessus (a priori $F$ peut être vide, on a alors $k=0$).

Par hypothèse, $d_1+ \cdots +d_k <n$. En particulier, $d_i<n$ et $f^{(d_i)}(t_i)\neq 0$ pour tout $i$. Pour~$t \in I$, posons
\[
P_i(t)=
\begin{dcases}
(t-t_i)^{d_i} & \text{ si } t_i <T\\
(t_i-t)^{d_i} & \text{ si } t_i =T.
\end{dcases}
\]

Les polynômes $P_i$ prennent des valeurs positives ou nulles sur $I$. La continuité de $f$ et la formule de Taylor entraînent alors l'existence de nombres $\delta >0$ et $c>0$ tels que
\begin{equation}\label{200119a}
\abs{f(t)} \ioe 
\begin{dcases}
1-cP_i(t) & (1\ioe i \ioe k\; ; \; t\in I, \, \abs{t-t_i} < \delta)\\
1-c & (t\in I, \,\min_{1\ioe i \ioe k} \abs{t-t_i} \soe \delta).
\end{dcases}
\end{equation}

Soit $P=P_1\cdots P_k$ et $K$ un majorant commun à toutes les fonctions polynomiales $P/P_i$, pour~$1\ioe i \ioe k$, et à $P$ elle-même, sur le segment $I$ (si $F$ est vide et $k=0$, on a $P=1$).

D'après \eqref{200119a}, pour $t \in I$, on a, d'une part, si $1 \ioe i \ioe k$ et $\abs{t-t_i} < \delta$, 
\[
\abs{f(t)} \ioe 1-cP_i(t) \ioe 1-\frac cK P(t),
\]
et, d'autre part, si $\min_{1\ioe i \ioe k} \abs{t-t_i} \soe \delta$,
\[
\abs{f(t)} \ioe 1-c \ioe 1-\frac cK P(t).
\]

On a donc $\abs{f}\ioe 1-\eps P$, avec $\eps=c/K$. Comme $P$ est une fonction polynomiale de degré~$<n$, sa dérivée $n$\up{e} est nulle. On en déduit que les deux fonctions $f \pm \eps P$ appartiennent à $\Lcal_n$, ce qui contredit le fait que $f$ est un point extrême de $\Lcal_n$. Cette contradiction démontre $(i)$.

\smallskip

Démontrons maintenant que $f$, point extrême de $\Lcal_n$, vérifie la condition $(ii)$. L'ensemble complémentaire $I \setminus F$ est un ouvert de $I$ ; c'est donc une réunion dénombrable d'intervalles~$I_j \subset I$. Pour montrer que
\[
 \mes\set{t \in I \setminus F, \, \abs{f^{(n)} (t)} \neq 1} =0,
\]
il suffit de montrer que, pour chaque $j$,
\[
\mes\set{t \in I_j, \, \abs{f^{(n)}(t)} \neq 1} =0,
\]
ou encore 
\begin{equation}\label{190223a}
\mes\set{t \in J, \, \abs{f^{(n)} (t)} \neq 1} =0,
\end{equation}
si $J$ est un segment quelconque inclus dans $I_j$. Observons tout de suite que, $\abs{f}$ étant continue et $<1$ sur le compact $J$, il existe $\eps >0$ tel que
\[
\abs{f(t)} \ioe 1-\eps \quad (t \in J).
\]

Comme $f \in \Lcal_n$, on a $\abs{f^{(n)} (t)}\ioe 1$ presque partout, donc~\eqref{190223a} équivaut à
\begin{equation}\label{190223b}
\forall \, \delta >0, \quad \mes\set{t \in J, \, \abs{f^{(n)}(t)} \ioe1-\delta} =0.
\end{equation}

Fixons $\delta >0$, et posons
\[
K=\set{t \in J, \, \abs{f^{(n)} (t)} \ioe1-\delta} .
\]
C'est un ensemble mesurable, et il s'agit de montrer que $\mes(K)=0$. Nous allons raisonner par l'absurde, en supposant que $\mes (K) >0$. Notons alors $w$ la fonction indicatrice de $K$,
\[
w(t)=[t \in K]  \quad (0 \ioe t \ioe T).
\]

On a $w \soe 0$ et $\int_0^T w(t) \, dt =\mes(K) >0$. La théorie des polynômes orthogonaux affirme alors l'existence d'une suite de polynômes $(p_k)_{k \in \Nat}$, à coefficients réels, où $p_k$ est de degré $k$, tels que
\begin{equation}\label{190227a}
\int_0^T p_j(t)\,p_k(t) \, w(t) \, dt =[j=k] \quad (j, k \in \Nat)
\end{equation}
(cf. \cite{zbMATH03477793}, \S 2.2, pp. 25-26).

Pour $t \in I$, posons
\begin{align*}
g_0(t)&=\int_0^t p_n(u)\, w(u) \,du \\
g_1(t)&=\int_0^t g_0(u) \,du=\int_0^t p_n(u)\, w(u)\,(t-u)\, du\\
&\vdots\\
h(t)=g_{n-1}(t)&=\int_0^t g_{n-2}(u) \,du=\int_0^t p_n(u)\, w(u)\,\frac{(t-u)^{n-1}}{(n-1)!}\, du.
\end{align*}

La fonction $h=g_{n-1}$ est de classe $\Ccal^{n-1}$ et $h^{(n-1)}=g_0$. Comme $g_0$ est l'intégrale de la fonction mesurable et bornée $p_nw$, elle est lipschitzienne sur~$I$, de rapport $\nu=\norm{g'_0}_{\infty}=\norm{p_nw}_{\infty}$. Par conséquent, $h^{(n-1)}$ est dérivable presque partout sur $I$, avec $ h^{(n)}=g'_0=p_nw$ presque partout. Observons que $ h^{(n)}$ n'est pas égale presque partout à la fonction nulle, et notons $\xi=\norm{h}_{\infty}$.

\smallskip

Soit $\alpha =\inf J$ et $\beta =\sup J$. Pour $u < \alpha$, on a $w(u)=0$ donc la fonction $h$ est nulle sur $[0,\alpha]$. Pour $u> \beta$, on a aussi $w(u)=0$ donc, pour $t \in [\beta,T]$,
\[
h(t)=\int_0^t p_n(u)\, w(u)\,\frac{(t-u)^{n-1}}{(n-1)!}\, du=\int_0^T p_n(u)\, w(u)\,\frac{(t-u)^{n-1}}{(n-1)!}\, du=0,
\]
d'après \eqref{190227a}, en écrivant, pour $t$ fixé, le polynôme $(t-u)^{n-1}/(n-1)!$, en la variable $u$, comme combinaison linéaire des polynômes $p_k(u)$ pour $0 \ioe k <n$.

Ainsi, la fonction $h$ s'annule sur $I\setminus J$, et
\begin{align*}
h^{(n)}(t)&=p_n(t)w(t)=0 \quad \text{pour presque tout } t \in I\setminus K\, ;\\
\abs{h^{(n)}(t)}&\ioe \nu \quad \text{pour presque tout } t \in I.
\end{align*}

Soit $ \eta >0$. On a
\begin{align*}
f(t) \pm \eta h(t) &=f(t) \quad ( t \in I\setminus J)\, ;\\
\abs{f(t) \pm \eta h(t)} &\ioe 1-\eps +\eta \,\xi \quad ( t \in J)\, ;\\
f^{(n)}(t) \pm \eta h^{(n)}(t) &=f^{(n)}(t) \quad \text{pour presque tout } t \in I\setminus K\, ;\\
\abs{f^{(n)}(t) \pm \eta h^{(n)}(t)} &\ioe 1-\delta +\eta\, \nu \quad \text{pour presque tout } t \in K.
\end{align*}

Par conséquent, si $\eta=\min(\delta/\nu,\eps/\xi)$, les fonctions $f \pm \eta h$ appartiennent à $\Lcal_n$. Comme $h$ n'est pas la fonction nulle, cela contredit le fait que $f$ est un point extrême de $\Lcal_n$. Cette contradiction prouve que $\mes (K)=0$. Par suite, $f$ vérifie $(ii)$.

\medskip

Passons à la réciproque. En supposant que $f\in \Lcal_n$ vérifie $(i)$, $(ii)$, et 
\begin{equation}\label{200114a}
2f=g+h,
\end{equation}
où $g,\, h \in \Lcal_n$, il nous faut montrer que $g=h$.

Nous allons d'abord montrer que les hypothèses $(ii)$ et \eqref{200114a}, et l'appartenance de $f$, $g$ et $h$ à $\Lcal_n$, entraînent que $g^{(n)}=h^{(n)}$ presque partout, \cad que $g-h$ est une fonction polynomiale, de degré $ <n$. Notons
\[
E=\set{t_0\in I, \, \text{$g$ ou $h$ n'est pas $n$ fois dérivable en $t_0$, ou bien } g^{(n)}(t_0)\neq h^{(n)}(t_0)}.
\]
Il s'agit de montrer que la mesure de $E$ est nulle. \'Ecrivons
\[
E=(E\cap F) \cup\big(E\cap(I\setminus F)\big).
\]
 
Montrons d'abord que $E\cap(I\setminus F)$ est de mesure nulle. Soit~$E_0$ l'ensemble des éléments~$t_0 \in I \setminus F$ vérifiant au moins l'une des conditions suivantes :

\smallskip

$\bullet$ $g^{(n)}(t_0)$ n'existe pas ;

$\bullet$ $h^{(n)}(t_0)$ n'existe pas ;

$\bullet$ $f^{(n)}(t_0)$ existe et $\abs{f^{(n)}(t_0)}\neq 1$ ;

$\bullet$ $g^{(n)}(t_0)$ et $h^{(n)}(t_0)$ existent, et $\max(\abs{g^{(n)}(t_0)},\abs{h^{(n)}(t_0)}) >1$.

\smallskip

D'après $(ii)$ et le fait que $g$ et $h$ appartiennent à $\Lcal_n$, l'ensemble $E_0$ est de mesure nulle. Si l'élément $t_0$ de $I\setminus F$ n'appartient pas à $E_0$, alors $g^{(n)}(t_0)$ et $h^{(n)}(t_0)$ existent, donc aussi $f^{(n)}(t_0)$, d'après \eqref{200114a}, et
\[
1=\lvert f^{(n)}(t_0)\rvert =\demi\lvert g^{(n)}(t_0)+h^{(n)}(t_0) \rvert  \ioe \demi\big( \lvert g^{(n)}(t_0)\rvert +\lvert h^{(n)}(t_0) \rvert \big)\ioe 1,
\]
d'où $g^{(n)}(t_0) =h^{(n)}(t_0)$, autrement dit $t_0 \not \in E$. On a donc~$E\cap(I\setminus F)\subset E_0$, et la mesure de~$E\cap(I\setminus F)$ est nulle.

Pour montrer que la mesure de $E\cap F$ est nulle, écrivons $F=F' \cup (F\setminus F')$, où $F'$ est l'ensemble des points d'accumulation de $F$. D'une part, l'ensemble $F\setminus F'$ est constitué de points isolés, donc est au plus dénombrable. D'autre part, si $t_0 \in F$, l'égalité~$\abs{g(t_0) +h(t_0)}=2$ et les inégalités~$\abs{g(t_0)}\ioe 1$ et $\abs{h(t_0)}\ioe 1$ entraînent l'égalité $g(t_0)=h(t_0)$. Si, de plus, $t_0 \in F'$, et si~$g^{(n)}(t_0)$ et $h^{(n)}(t_0)$ existent, la formule de Taylor en $t_0$ entraîne que $(g-h)^{(n)}(t_0)=0$ ; la mesure de $E\cap F'$ est donc aussi nulle. 

Nous avons bien montré que $g-h$ est polynomiale, de degré $<n$.

\smallskip

Pour en conclure que $g=h$, en supposant aussi $(i)$ vérifiée, il nous suffit de montrer que, pour tout $t_0 \in F$, les égalités
\begin{equation}\label{200121a}
f'(t_0)=\dots=f^{(k)}(t_0)=0
\end{equation}
pour un certain $k<n$, entraînent les égalités
\begin{equation*}
(g-h)(t_0)=(g-h)'(t_0)=\dots=(g-h)^{(k)}(t_0)=0.
\end{equation*}
En effet, l'hypothèse $(i)$ entraînera alors que la fonction polynomiale $g-h$, de degré $<n$, a, compte tenu des multiplicités, au moins $n$ zéros, et est donc nécessairement la fonction nulle.

D'une part, on a $g(t_0)=h(t_0)$, puisque $t_0 \in F$. D'autre part, nous allons montrer que \eqref{200121a} entraîne en fait que
\begin{equation}\label{200121c}
g'(t_0)=h'(t_0)=\dots=g^{(k)}(t_0)=h^{(k)}(t_0)=0.
\end{equation}

En effet, si \eqref{200121c} est fausse, soit $j$ le plus petit nombre entier positif tel que $g^{(j)}(t_0)\neq 0$ ou~$h^{(j)}(t_0)\neq 0$. On a donc $1 \ioe j \ioe k$. Pour fixer les idées, posons $c=g^{(j)}(t_0)/j!\neq 0$. La formule de Taylor nous donne alors
\begin{align}
g(t) &=g(t_0) +c(t-t_0)^j +o\big((t-t_0)^j\big) \quad (t\vers t_0, \, t \in I)\label{200121d}\\
f(t) &= f(t_0)+o\big((t-t_0)^j\big) \quad (t\vers t_0, \, t \in I)\label{200121e}\expli{d'après \eqref{200121a}}\\
h(t) &=h(t_0) -c(t-t_0)^j +o\big((t-t_0)^j\big) \quad (t\vers t_0, \, t \in I)\label{200121f},
\end{align}
où \eqref{200121f} résulte de \eqref{200121d}, \eqref{200121e}, et de la relation $h=2f-g$. Comme $c\neq 0$ et $g(t_0)=h(t_0)=\pm 1$, les relations \eqref{200121d} et \eqref{200121f} entraînent que $\max\big(\abs{g(t)},\abs{h(t)}\big) >1$ dans un certain voisinage (épointé) de $t_0$ dans $I$. Cela contredit l'appartenance de $g$ et $h$ à $\Lcal_n$, et cette contradiction achève la démonstration de la proposition.\fin

\smallskip

Dans la suite, nous noterons $\Fcal_n$ l'ensemble des points extrêmes de $\Lcal_n$, caractérisés dans la proposition \ref{200118a}. 

\smallskip

Notons que, dans le cas $n=2$, la condition $(i)$ équivaut à la conjonction des deux conditions suivantes :

$\bullet$ le fermé $F$ est non vide ;

$\bullet$ si $F=\set{0}$ (resp. $F=\set{T}$), on a $f'(0)=0$ (resp. $f'(T)=0$).

\smallskip

Pour ce cas $n=2$, donnons un exemple de fonction, point extrême de $\Lcal_2$, que nous utiliserons au \S\ref{190605c}.

\begin{prop}\label{190531a}
Soit $t_0$ et $T>0$, tels que $0\ioe t_0 \ioe T$. On pose
\[
C=\frac 2T+\frac{t_0^2+(T-t_0)^2}{2T}\virg
\]
et on suppose que 
\begin{equation}\label{190531b}
C \soe \max(t_0,T-t_0).
\end{equation}

Alors la fonction $f$ définie sur $[0,T]$ par 
\begin{equation}\label{190531c}
f(t)=
\begin{dcases}
-1+(C-t_0)t+t^2/2 & ( 0 \ioe t \ioe t_0)\\
1+(C-T+t_0)(t-T)-(t-T)^2/2 & ( t_0 \ioe t \ioe T),
\end{dcases}
\end{equation}
est un point extrême de $\Lcal_2$ tel que $f'(t_0)=C$.
\end{prop}
\dem

La constante $C$ est choisie pour que les deux formules donnent la même valeur pour $f(t_0)$, et on a $f''(t)=\pm 1$ si $t \neq t_0$. Ensuite, les dérivées à droite et à gauche de $f$ en~$t=t_0$ ont pour valeur commune $C$ (si $t_0=0$ ou $T$, il n'y a pas besoin de raccordement). Enfin, la condition~\eqref{190531b} entraîne que $f$ est strictement croissante sur $[0,T]$, avec~$f(0)=-1$ et~$f(T)=1$ ; on a donc~$\abs{f} \ioe 1$ sur~$[0,T]$. Ainsi, $f \in \Lcal_2$, l'ensemble $F$ de la proposition \ref{200118a} est égal à $\set{0,T}$ et les points $(i)$ et~$(ii)$ de cette proposition sont vérifiés.\fin

\section{Le noyau de Peano}\label{190322b}

Nous montrons dans ce paragraphe que les dérivées intermédiaires $f^{(k)}$ des fonctions $f \in~\!\Lcal_n(T)$ sont bornées, en valeur absolue, par une quantité dépendant uniquement de $k$, $n$ et $T$. 

Commençons par le cas $n=2$. La majoration de la valeur absolue de la dérivée d'une fonction de $\Lcal_2(T)$ s'appuie sur l'identité suivante, déjà implicite dans la démonstration du Satz 1, p. 45, de~\cite{zbMATH02621634}.

\begin{prop}\label{190321a}
Soit $f : [0,T] \vers \Real$ une fonction dérivable, dont la dérivée est absolument continue sur $[0,T]$. Pour $0 \ioe x \ioe T$, on a
\[
f'(x)=\frac{f(T)-f(0)}{T}+\int_{0}^T K(x,t\, ; \, T) \, f''(t) \, dt,
\]
où, pour $0 \ioe x ,t \ioe T$, $x \neq t$,
\[
K(x,t\, ; \, T)=
\begin{dcases}
t/T &(0 \ioe t < x)\\
(t-T)/T & (x < t \ioe T).
\end{dcases}
\]
\end{prop}
\dem

Observons que nous ne supposons pas que $f''$ (définie presque partout) est essentiellement bornée. Mais, par l'hypothèse d'absolue continuité de $f'$, cette dérivée seconde est, en tout cas, intégrable au sens de Lebesgue sur $[0,T]$ et les intégrations par parties qui vont suivre sont bien licites.

Ainsi, on voit que
\begin{align*}
\int_{0}^x ( -t) \, f''(t) \, dt &=f(x)-f(0)+ -xf'(x)   \\
\int_x^{T} (T -t) \, f''(t) \, dt &= f(T)-f(x) -f'(x) (T-x).
\end{align*}

En ajoutant ces deux égalités, on obtient l'énoncé.\fin

\smallskip

La norme $L^1$ (par rapport à la variable $t$) du noyau $K(x,t\, ; \, T)$ vaut
\begin{align}
\int_{0}^{T} \abs{K(x,t\, ; \, T)} \, dt &=\int_{0}^x \frac{t}{T}\, dt+\int_x^{T}\frac{T-t}{T}\, dt\notag\\
&=\frac{x^2+(T-x)^2}{2T}\cdotp\label{190529a}
\end{align}
Elle est maximale pour $x=0$ et $x=T$, prenant alors la valeur $T/2$.

\begin{prop}\label{190321b}
Pour toute fonction $f \in \Lcal_2(T)$, on a $\norm{f'}_{\infty} \ioe T/2+2/T$.
\end{prop}
\dem

Soit $f \in \Lcal_2(T)$ et $x \in [0,T]$. D'après la proposition \ref{190321a}, on a
\begin{align}
\abs{f'(x)} & \ioe \abs{f(T)-f(0)}/T +\int_0^T \abs{K(x,t\, ; \, T)} \cdot \abs{f''(t)} \, dt\label{200203a}\\
&\ioe 2/T+ \int_0^T \abs{K(x,t\, ; \, T)} \, dt\notag\\
&=2/T +(x^2+(T-x)^2)/2T\expli{d'après \eqref{190529a}}\notag\\
&= 2/T+T/2-x(T-x)/T \label{190605a}\\
&\ioe 2/T +T/2.\fine\notag
\end{align}

\smallskip

Nous déterminerons, à la proposition \ref{190605b} du \S\ref{190605c}, les cas où la majoration \eqref{190605a} est optimale, et nous verrons, à la proposition \ref{190605d} du \S\ref{190216a}, que la majoration $\norm{f'}_{\infty} \ioe T/2+2/T$ est optimale si $T\ioe 2$.

\smallskip

Dans l'article \cite{zbMATH02621356}, publié la même année, 1913, où Landau présentait son travail à la London Mathematical Society, Peano donna une formule générale dont la proposition \ref{190321a} ci-dessus est un cas particulier. Il y considérait les intégrales successives de la fonction de Heaviside, définies sur~$\Real$ par
\begin{align*}
H_0(x) &= [x \soe 0]\\
H_m(x)&=\int_0^x H_{m-1}(t) \, dt=[x \soe 0] \cdot x^m/m! \quad ( m \soe 1).
\end{align*}

La fonction $H_m$ est de classe $\Ccal^{\infty}$ sur $\Real \setminus\set{0}$, et, si $m \soe 1$, de classe $\Ccal^{m-1}$ sur $\Real$. On a~$H_m^{(j)}=H_{m-j}$ pour $0 \ioe j \ioe m$.

Voici un énoncé du théorème de Peano, qui n'est pas le plus général possible, mais qui nous suffira dans cet article. 
\begin{prop}\label{200128a}
Soit $T>0$, $k \in \Nat$ et, pour $i=1,\dots,k$, 
\[
\alpha_i \in [0,T] \quad ; \quad  \lambda_i \in \Real \quad ; \quad m_i \in \set{0,\dots, n-1}. 
\]

Considérons la forme linéaire $L$ définie par
\[
L(f)=\sum_{1\ioe i \ioe k} \lambda_if^{(m_i)}(\alpha_i),
\]
pour toute fonction $f : [0,T] \vers \Real$ telle que les nombres dérivés figurant dans cette somme existent (en particulier, c'est le cas si $f \in V_n([0,T])$). On suppose que $L(f)=0$ si $f$ est une fonction polynomiale, de degré $<n$. 

On a alors, pour toute fonction $f \in V_n([0,T])$,
\begin{equation}\label{200124a}
L(f)=\int_0^T K(t) f^{(n)}(t) \, dt,
\end{equation}
où la fonction $K$, appelée \emph{noyau de Peano} associé à la forme linéaire $L$ pour l'ordre $n$, est définie par la formule
\[
K(t)=L(H_{n-1,t}),
\]
où $H_{n-1,t}$ est la fonction définie sur $\Real$ par $H_{n-1,t}(x)=H_{n-1}(x-t)$.
 
Le noyau de Peano $K(t)$ est défini pour $t \in [0,T]\setminus\set{\alpha_i, \;m_i=n-1}$. C'est une fonction polynomiale par morceaux, dont la valeur absolue est bornée sur son ensemble de définition par la quantité
\[
\sum_{1 \ioe i \ioe k} \abs{\lambda_i }\frac{T^{n-1-m_i}}{(n-1-m_i)!}\cdotp
\]
\end{prop}
\dem

Soit $f \in V_n([0,T])$ et
\[
P(t)=f(0)+f'(0)t + \cdots + f^{(n-1)}(0)t^{n-1}/(n-1)!
\]
le polynôme de Taylor d'ordre $n-1$ de $f$ en $0$. Posons $g=f-P$, de sorte que $g \in V_n([0,T])$. On a l'égalité~$f^{(n)}=g^{(n)}$ presque partout, et $L(f)=L(g)$.

Appliquons à chaque terme de $L(g)$ la formule de Taylor avec reste intégral. Comme les dérivées successives de $g$ au point $0$ sont nulles jusqu'à l'ordre $n-1$, on a
\begin{equation}\label{200124b}
g^{(m_i)}(\alpha_i)=\int_0^{\alpha_i} \frac{(\alpha_i-t)^{n-1-m_i}}{(n-1-m_i)!}\,g^{(n)}(t) \, dt=\int_0^T H_{n-1,t}^{(m_i)}(\alpha_i)\, f^{(n)}(t) \, dt
\end{equation}

On obtient \eqref{200124a} en multipliant \eqref{200124b} par $\lambda_i$ et en sommant pour $1 \ioe i \ioe k$. La dernière assertion résulte de l'expression
\begin{equation*}
K(t) =\sum_{1 \ioe i \ioe k} \lambda_i \frac{(\alpha_i-t)^{n-1-m_i}}{(n-1-m_i)!} \cdot[t \ioe \alpha_i].\fine
\end{equation*}

\smallskip

Ainsi, la proposition \ref{190321a} explicite le noyau de Peano associé à la forme linéaire
\[
f \mapsto f'(x)-\big(f(T)-f(0)\big)/T
\]
pour l'ordre $2$. Un exposé détaillé de la théorie du noyau de Peano est donné par Sard au chapter~1 de \cite{zbMATH03187028} (cf. en particulier le \S 43, p. 25).

\smallskip

La proposition \ref{190321b} peut être généralisée, au moins qualitativement, grâce au noyau de Peano. Pour cela, il est commode de revenir à une formulation pour l'ensemble $\Lcal_n(a,b \, ; [0,T])$ au lieu de $\Lcal_n(T)$.

\begin{prop}\label{200402a}

Soit $k$ et $n$ des nombres entiers, tels que $n \soe 2$ et $0 <k<n$. Il existe deux constantes positives $A_{n,k}$ et $B_{n,k}$, telles que, pour $T>0$, $a >0$, $b>0$, et $f \in \Lcal_n(a,b \, ; [0,T])$, on ait
\begin{equation}\label{200402b}
\|f^{(k)}\|_{\infty} \ioe A_{n,k}\,aT^{-k} +B_{n,k}\,bT^{n-k}.
\end{equation}
\end{prop}
\dem

L'argumentation qui suit est une variante de celle donnée par deVore et Lorentz dans \cite{zbMATH00477682}, Theorem~5.6, p. 38. 

Tout d'abord, l'existence de constantes $A_{n,k}$ et $B_{n,k}$ telles que \eqref{200402b} soit vraie pour $T=1$ (et tous $a$, $b$, $f$), entraîne \eqref{200402b} pour tout $T$. En effet, si $f \in \Lcal_n(a,b \, ; [0,T])$, alors la fonction $g$, définie sur $[0,1]$ par $g(u)=f(Tu)$, appartient à $\Lcal_n(a,bT^{n} \, ; [0,1])$, donc \eqref{200402b} avec $T=1$ entraîne
\[
T^k\|f^{(k)}\|_{\infty}=\|g^{(k)}\|_{\infty} \ioe A_{n,k}\,a +B_{n,k}\,bT^{n},
\]
d'où \eqref{200402b}.

\smallskip

Ensuite, pour démontrer \eqref{200402b} quand $T=1$, commençons en fixant arbitrairement $n$ éléments de $[0,1]$,
\[
0 \ioe \alpha_1 < \dots < \alpha_n \ioe 1,
\]
par exemple $\alpha_i=i/n$.

Si $0 \ioe x\ioe 1$, déterminons les coefficients $\lambda_i$, pour $1 \ioe i \ioe n$, de sorte que la forme linéaire
\[
L_x(f)=f^{(k)}(x)-\sum_{1\ioe i \ioe n} \lambda_if(\alpha_i)
\]
s'annule lorsque $f$ est une fonction polynomiale, de degré $<n$. Les $\lambda_i$ sont solutions d'un système linéaire $n \times n$, dont le déterminant est le déterminant de Vandermonde formé avec les $\alpha_i$ ; les formules de Cramer montrent que les valeurs absolues de ces solutions sont bornées par des quantités, dépendant de $k$ et $n$, mais indépendantes de $x \in [0,1]$.

D'après la proposition \ref{200128a}, le noyau de Peano $K_x$, associé à la forme linéaire $L_x$ pour l'ordre $n$ a également une valeur absolue bornée par une quantité, disons $B_{n,k}$, dépendant de $k$ et $n$, mais indépendante de $x \in [0,1]$.

Comme
\begin{align*}
f^{(k)}(x)&=\sum_{1\ioe i \ioe n} \lambda_if(\alpha_i) +L_x(f)\\
&=\sum_{1\ioe i \ioe n} \lambda_if(\alpha_i) +\int_0^1 K_x(t) f^{(n)}(t) \, dt,
\end{align*}
on en déduit que
\[
\lvert f^{(k)}(x)\rvert \ioe A_{n,k} a + B_{n,k} b,
\]
si $f \in \Lcal_n(a,b \, ; [0,1])$, où
\begin{equation*}
A_{n,k}=\sum_{1\ioe i \ioe n} \lvert\lambda_i\rvert. \fine
\end{equation*}

\smallskip

Posons $A_{n,k}^*=A_{n,k}/(n-k)$ et $B_{n,k}^*=B_{n,k}/k$. La valeur de $T$ pour laquelle le majorant dans~\eqref{200402b} est minimal est $T_0=( A_{n,k}^*a/B_{n,k}^*b)^{1/n}$. Pour $T\soe T_0$, tout élément $t$ de $[0,T]$ appartient à un segment de longueur $T_0$, inclus dans $[0,T]$, donc on peut majorer $\lvert f^{(k)}(t)\rvert$ en remplaçant $T$ par $T_0$ dans \eqref{200402b}. On obtient ainsi
\begin{equation*}\
\|f^{(k)}\|_{\infty} \ioe C_{n,k}^*\,a^{1-k/n}\,b^{k/n} \quad \big(f \in \Lcal_n(a,b\, ; [0,T]), \; 0<k<n \big),
\end{equation*}
où $C_{n,k}^*=nA_{n,k}^{*(1-k/n)}B_{n,k}^{*k/n}$.

\smallskip

Des expressions explicites de constantes $A_{n,k}$ et $B_{n,k}$, admissibles dans \eqref{200402b}, sont indiquées au \S\ref{200202a} ci-dessous.

\section{Compacité}\label{190322d}

En notant encore $I=[0,T]$, considérons l'espace vectoriel~$W=W_n(I)$ constitué des fonctions~$f\in V_n(I)$  telles que $f^{(n-1)}$ est lipschitzienne, autrement dit, telles que $\|f^{(n)}\|_{\infty} <\infty$. C'est un cas particulier d'espace de Sobolev. Considérons deux normes sur $W$ :
  \begin{align*}
N_1(f)&=\max(\norm{f}_{\infty},\|f'\|_{\infty},\dots, \|f^{(n-1)}\|_{\infty})\\
N_2(f)&=\max(\norm{f}_{\infty},\|f'\|_{\infty},\dots, \|f^{(n-1)}|_{\infty},\|f^{(n)}\|_{\infty})
\end{align*}

La topologie de $(W,N_2)$ est plus fine que celle de~$(W,N_1)$. L'espace vectoriel normé $(W,N_2)$ est complet, l'espace $(W,N_1)$ ne l'est pas. 

\begin{prop}
L'ensemble $\Lcal_n=\Lcal_n(1,1 \, ; \, [0,T])$ est

$\bullet$ compact, donc fermé et borné, dans $(W,N_1)$ ;

$\bullet$ fermé et borné, mais non compact, dans $(W,N_2)$.
\end{prop}
\dem

D'une part, d'après la proposition \ref{200402a}, $\Lcal_n$ est borné dans $(W,N_1)$ et dans $(W,N_2)$. Comme les évaluations ponctuelles, $f \mapsto f(t_0)$ et $f \mapsto f^{(n-1)}(t_0)$, sont continues sur $(W,N_1)$ pour tout~$t_0 \in I$, l'ensemble $\Lcal_n$, qui est défini par les inégalités
\[
\abs{f(t)} \ioe 1 \quad (t \in I)  ; \quad \abs{f^{(n-1)}(t_1)-f^{(n-1)}(t_0)} \ioe \abs{t_1-t_0} \quad (t_0,\, t_1 \in I),
\]
est un fermé de $(W,N_1)$, et donc aussi de $(W,N_2)$. Dans ce dernier espace, $\Lcal_n$ contient un voisinage ouvert de~$0$, et n'est donc pas compact, puisque $W$ est de dimension infinie.

Munissons $\Real^n$ de la norme $\norm{(x_0,\dots,x_{n-1})}=\max(\abs{x_0},\dots,\abs{x_{n-1}})$, et considérons l'espace~$\Ccal$ des fonctions continues, définies sur $I$, à valeurs dans $\Real^n$. Muni de la norme uniforme, $\Ccal$ est un espace de Banach. L'application linéaire $f \mapsto(f,f',\dots,f^{(n-1)})$ est un isomorphisme d'espaces vectoriels normés entre $(W,N_1)$ et un sous-espace (non fermé) de $\Ccal$, que nous noterons $\Ccal_1$.

Soit $\Lcal'_n \subset \Ccal_1$ l'image de $\Lcal_n$ par cet isomorphisme. Observons d'abord que $\Lcal'_n$ est un fermé, non seulement de $\Ccal_1$, mais aussi  de~$\Ccal$. En effet, si la suite $(f_k)_{k \in \Nat}$ d'éléments de $\Lcal_n$ est telle que
\[
(f_k,f'_k,\dots,f^{(n-1)}_k) \vers (g_0,g_1,\dots,g_{n-1}) \in \Ccal \quad (k \vers \infty),
\]
alors, puisque la convergence est uniforme, la fonction $g_0$ est $(n-1)$ fois dérivable, et l'on a~$g_0^{(i)}=g_i$ pour $0 \ioe i \ioe n-1$. De plus, la fonction $g^{(n-1)}$, limite uniforme de la suite de fonctions lipschitziennes de rapport $1$, $(f^{(n-1)}_k)$, est elle-même lipschitzienne de rapport $1$.

D'après la proposition \ref{200402a}, les normes des fonctions appartenant à $\Lcal'_n$ sont uniformément bornées par une quantité $\lambda=\lambda_n(T) \in [1,\infty[$, et toutes ces fonctions sont lipschitziennes, de rapport~$\lambda$. L'ensemble~$\Lcal'_n$ est donc une partie équicontinue de $\Ccal$ et, pour tout $t_0 \in I$, l'ensemble~$\set{G(t_0), \, G \in \Lcal'_n}$ est borné dans $\Real^n$, donc a une adhérence compacte.

D'après le théorème d'Ascoli (cf., par exemple, \cite{zbMATH03703614}, \S 6.3, p. 87-88), $\Lcal'_n$ a une adhérence compacte dans $\Ccal$. Comme $\Lcal'_n$ est fermé, c'est un compact, et son image isomorphe $\Lcal_n$ dans l'espace~$(W,N_1)$ est aussi compacte.\fin

\smallskip

Le théorème de Krein et Milman (cf. \cite{zbMATH01022519}, Theorem 3.23, p. 75), entraîne donc la proposition suivante.
\begin{prop}
L'ensemble $\Lcal_n$ est l'enveloppe convexe fermée dans $(W,N_1)$ de l'ensemble $\Fcal_n$ de ses points extrêmes. 
\end{prop}

Cela étant, le Theorem 3.2, p. 1159 de \cite{zbMATH03256144} incite à se demander si une propriété plus forte est vraie ou non.
\begin{ques}
L'ensemble $\Lcal_n$ est-il l'enveloppe convexe fermée de $\Fcal_n$ dans $(W,N_2)$?
\end{ques}

\section{Problèmes extrémaux dans $\Lcal_n$ : généralités}\label{190322e}

\subsection{Fonctions extrémales}

Plusieurs variantes du problème de Landau étudiées dans la littérature sont du type suivant. On se donne une fonction $\Phi$, à valeurs réelles, définie sur $\Lcal_n=\Lcal_n(T) $, convexe, et continue pour la topologie de $(W,N_1)$, définie au paragraphe précédent. Il s'agit de déterminer la borne supérieure $\sup \Phi$, et l'ensemble
\[
\Kcal(\Phi)=\set{f \in \Lcal_n, \, \Phi(f)=\sup \Phi}.
\]

Les éléments de $\Kcal(\Phi)$ sont les \emph{fonctions extrémales} pour $\Phi$.

Avant d'énoncer la proposition suivante, qui rassemble quelques faits classiques, rappelons la définition de partie extrême. Si $A\subset B\subset E$, où $E$ est un espace vectoriel réel, on dit que $A$ est extrême dans $B$ si
\[
\forall\, b_1,b_2\in B, \; \forall t \in \, ]0,1[, \; tb_1+(1-t)b_2 \in A \implique b_1 \in A \text{ et } b_2 \in A.
\]

Ainsi, par exemple, les points extrêmes de $B$ correspondent aux parties extrêmes de $B$ ne contenant qu'un élément. La relation binaire {\og être extrême dans\fg} est une relation d'ordre sur l'ensemble $\Pcal(E)$ des parties de $E$. En particulier, si $A$ est extrême dans $B$, et si $C$ désigne l'ensemble des points extrêmes de~$B$, alors l'ensemble des points extrêmes de $A$ est $A \cap C$.

\begin{prop}\label{190615a}
L'ensemble $\Kcal=\Kcal(\Phi)$ est

$\bullet$ non vide ;

$\bullet$ compact pour la topologie de $(W,N_1)$ ;

$\bullet$ extrême dans $\Lcal_n$ ;

$\bullet$ contenu dans l'enveloppe convexe fermée de $\Fcal_n\cap \Kcal$ dans $(W,N_1)$.

\smallskip

Si, en outre, $\Phi$ est la restriction à $\Lcal_n$ d'une forme linéaire continue sur $(W,N_1)$, l'ensemble~$\Kcal$ est

$\bullet$ convexe, égal à l'enveloppe convexe fermée de $\Fcal_n\cap \Kcal$ dans $(W,N_1)$.

\end{prop}
\dem

Les deux premières assertions découlent de la continuité de $\Phi$ et de la compacité de $\Lcal_n$ pour la topologie de $(W,N_1)$. La troisième assertion découle de la convexité de $\Phi$.

Comme $\Kcal$ est une partie extrême de $\Lcal_n$, l'ensemble des points extrêmes de $\Kcal$ est $\Fcal_n\cap \Kcal$. La quatrième assertion résulte alors de la variante du théorème de Krein et Milman, énoncée comme Theorem 3.24, p. 76 de \cite{zbMATH01022519}.

Enfin , la cinquième assertion résulte de l'hypothèse supplémentaire de linéarité, et du théorème de Krein et Milman.\fin

\smallskip

Il résulte de la première et de la quatrième assertion de la proposition \ref{190615a} que $\Fcal_n\cap \Kcal$ est non vide : il existe un point extrême de $\Lcal_n$ où la borne supérieure de $\Phi$ est atteinte. Autrement dit, pour déterminer cette borne supérieure, on peut se contenter d'examiner la restriction de $\Phi$ à l'ensemble $\Fcal_n$ des points extrêmes de $\Lcal_n$, dont la description est l'objet de la proposition \ref{200118a}.

\subsection{Monotonie}

Une autre propriété simple, mais essentielle, est la décroissance des bornes supérieures \eqref{200425b} lorsque l'intervalle $I$ grandit. 

\begin{prop}\label{200425a}
Soit $a$ et $b$ des nombres réels strictement positifs, $k$ et $n$ des nombres entiers tels que $0<k<n$. Si $I$ est un intervalle fermé de $\Real$, posons
\[
\sigma_{\infty}(a,b \, ;\, k,n\,; \, I)=\sup_{f \in \Lcal_n (a,b\, ; \, I)} \|f^{(k)}\|_{\infty}.
\]

Si $J$ est un intervalle fermé de $\Real$ tel que $I \subset J$, on a 
\[
\sigma_{\infty}(a,b\,  ;\, k,n\,; \, J) \ioe \sigma_{\infty}(a,b \, ;\, k,n\,; \, I).
\]
\end{prop}
\dem

Soit $f\in \Lcal_n (a,b\, ; \, J)$ et $t_0 \in J$. Il existe un nombre réel $t_1$ tel que 
\[
t_0 \in t_1+I \subset J
\]
(on peut prendre $t_1=\pm\dist(t_0,I)$). La fonction $g : I \vers \Real$, définie par $g(t)=f(t+t_1)$, appartient à $\Lcal_n (a,b\, ; \, I)$, donc
\[
\lvert f^{(k)}(t_0)\rvert=\lvert g^{(k)}(t_0-t_1)\rvert \ioe \sigma_{\infty}(a,b \, ;\, k,n\, ; \, I).
\]

Comme $t_0$ est un élément quelconque de $J$, on a bien l'inégalité annoncée.\fin

\section{Le cas $n=2$}\label{200416a}

\subsection{Le problème extrémal de Landau pour $I=[0,T]$}\label{190216a}

Soit $T>0$. Pour simplifier l'écriture dans les cas où une seule valeur de $T$ intervient, nous posons $\Lcal=\Lcal_2(T)$.

Soit
\[
\sigma_{\infty}(T)=\sup_{f \in \Lcal}\, \norm{f'}_{\infty}.
\]

Des fonctions extrémales pour ce problème étant faciles à déterminer, les propositions \ref{190321b} et~\ref{200425a} vont nous fournir la valeur de $\sigma_{\infty}(T)$.

\begin{prop}\label{190605d}
 On a 
 \[
\sigma_{\infty}(T)=
\begin{dcases}
2/T+T/2 &(T \ioe 2) \\
2 & (T\soe 2).
\end{dcases}
\]
\end{prop}
\dem

Par la proposition \ref{190321b}, on a $\sigma_{\infty}(T) \ioe 2/T+T/2$. 

Si $T\ioe 2$, la fonction $f$ définie par
\begin{equation}\label{190620b}
f(t)=-t^2/2+(2/T+T/2)t-1
\end{equation}
appartient à $\Lcal$ et vérifie $f'(0)=2/T+T/2$, donc $\sigma_{\infty}(T)=2/T+T/2$.

\smallskip

Si $T > 2$, la proposition \ref{200425a} montre que $\sigma_{\infty}(T) \ioe \sigma_{\infty}(2)=2$. La fonction $f$ définie par
\begin{equation}\label{190621a}
f(t)=
\begin{dcases}
-t^2/2+2t-1 & (0 \ioe t \ioe 2)\\
1 & (t>2).
\end{dcases}
\end{equation}
appartient à $\Lcal$ et vérifie $f'(0)=2$, donc $\sigma_{\infty}(T)=2$.\fin

\subsection{Le problème extrémal de Landau pour $I=[0,\infty[$}\label{200423a}

La proposition suivante est un corollaire de la proposition \ref{190605d}.

\begin{prop}\label{190627g}
On a
\[
\sup_{f \in \Lcal_2(1,1 \, ;\, [0,\infty[\,)}\, \norm{f'}_{\infty} =2
\]
\end{prop}
\dem

La proposition \ref{200425a} montre que cette borne supérieure est $\ioe \sigma_{\infty}(2)=2$, borne atteinte si $f$ est définie par \eqref{190621a}.\fin

\smallskip

Schoenberg montre dans \cite{zbMATH03411359} (Theorem 9, p. 149) que toute fonction $f \in \Lcal_2(1,1 \, ;[0,\infty[\,)$, vérifiant $\norm{f'}_{\infty} =2$, est telle que
\[
f(t)=\pm (-t^2/2+2t-1) \quad (0 \ioe t \ioe 2).
\]

\subsection{Variante ponctuelle du problème extrémal de Landau}\label{190406a}

Fixons maintenant un point~$t_0 \in I=[0,T]$. Nous nous proposons de déterminer la quantité 
\[
\sigma(t_0,T)=\sup_{f \in \Lcal} f'(t_0),
\]
et une fonction $f$ extrémale pour ce problème, \cad telle que $f'(t_0)=\sigma(t_0,T)$. L'invariance de l'ensemble $\Lcal$ par l'application $f \mapsto -f$ montre que l'on a également
\[
\sigma(t_0,T)=\sup_{f \in \Lcal}\, \abs{f'(t_0)}.
\]

Ce problème est un cas particulier de celui étudié par Landau dans son article \cite{zbMATH02591405} de 1925. La condition sur $f''$ y est plus générale, étant de la forme $-2\alpha \ioe f'' \ioe 2\beta$, ou même simplement~$-2\alpha \ioe f'' $ ou $ f'' \ioe 2\beta$.

\subsubsection{Les fonctions $\fhi$ et $G$}

L'étude de ce problème s'appuie sur la considération des deux fonctions suivantes,
\begin{align*}
\fhi(x) &= \sqrt{2x^2+4} \, -x \quad ( x \soe 0)\\
G(x,y) &= \frac{2}{x+y}+\frac{x^2+ y^2}{2(x+y)} \quad ( x\soe 0, \, y \soe 0, \, x+y >0).
\end{align*}

Observons que $G$ est symétrique : $G(x,y)=G(y,x)$.

\begin{prop}\label{190521b}
La fonction $\fhi(x)$ est décroissante pour $0\ioe x \ioe \sqrt{2}$, croissante pour $x \soe \sqrt{2}$, avec $\fhi(0)=\fhi(4)=2$ et $\fhi(\sqrt{2})=\sqrt{2}$.

En particulier, on a $\fhi(x) \soe \sqrt{2}$ pour tout $x \soe 0$.
\end{prop}
\dem

Cela résulte du calcul de la dérivée,
\begin{equation*}
\fhi'(x) = -1+2/\sqrt{2+4/x^2}  \quad ( x > 0).\fine
\end{equation*}

\begin{prop}\label{190516a}
Pour $x ,y\soe 0$ tels que $x+y>0$, on a
\begin{align*}
G(x,y) & > y \quad (y < \fhi(x))\\
G(x,y) & = y \quad ( y = \fhi(x))\\
G(x,y) & < y \quad (y >\fhi(x)).
\end{align*}
\end{prop}
\dem

Cela résulte de l'identité
\begin{equation*}
\frac{y-G(x,y)}{x+y}=\frac 12-\frac{x^2+2}{(x+y)^2} \quad (x \soe 0, \, y \soe 0, \, x+y >0).\fine
\end{equation*}

\subsubsection{Détermination de $\sigma(t_0,T)$}\label{190605c}

Pour commencer, observons que l'application $f \mapsto g$, où $g(t)=-f(T-t)$, est une involution sur $\Lcal$, telle que $f'(t_0)=g'(T-t_0)$. On se ramène donc au cas $t_0\ioe T/2$, \cad $t_0 \ioe T-t_0$.

\begin{prop}\label{190603b}
Soit $f\in \Lcal$, $t_0 \in [0,T]$, et $h$, $h'$ tels que 
\[
h \soe 0, \; h' \soe 0, \; h+h'>0, \quad 0\ioe t_0-h \ioe t_0 \ioe t_0+h' \ioe T.
\]

Alors $f'(t_0) \ioe G(h,h').$
\end{prop}
\dem

En appliquant la proposition \ref{190321a} à la fonction $t \mapsto f(t_0-h+t)$ sur l'intervalle $[0,h+h']$, on obtient
\begin{align*}
f'(t_0) &=\frac{f(t_0+h')-f(t_0-h)}{h+h'}+\int_{0}^{h+h'}K(h,t\, ; \, h+h') \, f''(t_0-h+t) \, dt\\
&\ioe 2/(h+h')+\int_{0}^{h+h'}\abs{K(h,t\, ; \, h+h')}  \, dt \expli{car $f \in \Lcal$}\\
&=2/(h+h')+(h^2+h'^2)/2(h+h'),
\end{align*}
d'après \eqref{190529a}.\fin

\smallskip

Nous pouvons maintenant déterminer $\sigma(t_0,T)$. Pour chaque couple $(t_0,T)$, une fonction extrémale est mentionnée dans la démonstration.

\begin{prop}\label{190605b}
Soit $t_0 \in [0,T/2]$. On a 
\[
\sigma(t_0,T)=
\begin{dcases}
2/T+T/2-t_0(T-t_0)/T &(t_0 \ioe \sqrt{2} \,\text{ et }\, T \ioe \sqrt{2t_0^2+4}\,)\\
\sqrt{2t_0^2+4}\,-t_0 &(t_0 \ioe \sqrt{2} \,\text{ et }\, T > \sqrt{2t_0^2+4})\\
\sqrt{2}&(t_0 > \sqrt{2}).
\end{dcases}
\]
\end{prop}
\dem

D'après la proposition \ref{190603b}, on a toujours $\sigma(t_0,T) \ioe G(t_0,T-t_0)=2/T+T/2-t_0(T-t_0)/T$ (cf. \eqref{190605a}).

\smallskip

$\bullet$ Si $t_0 \ioe \sqrt{2} $ et $T-t_0 \ioe \fhi(t_0)$, on a aussi $t_0 \ioe \fhi(T-t_0)$, donc $G(t_0,T-t_0) \soe \max(t_0,T-t_0)$, d'après la proposition \ref{190516a}. La proposition \ref{190531a} donne alors explicitement une fonction $f_1\in \Lcal$ telle que $f_1'(t_0)=G(t_0,T-t_0)$.

\smallskip

$\bullet$ Si $t_0 \ioe \sqrt{2}$ et $T-t_0 > \fhi(t_0)$, le choix $h=t_0$, $h'=\fhi(t_0)$, dans la proposition \ref{190603b}, fournit l'inégalité
\[
\sigma(t_0,T) \ioe G\big(t_0,\fhi(t_0)\big)=\fhi(t_0).
\]

Pour montrer que $\sigma(t_0,T)=\fhi(t_0)$, considérons la fonction $f_2$ définie sur $[0,T]$ par 
\begin{equation}\label{200424a}
f_2(t)=
\begin{dcases}
-1+\big(\fhi(t_0)-t_0\big)t+t^2/2 & ( 0 \ioe t \ioe t_0)\\
1-(t-\fhi(t_0)-t_0)^2/2 & \big( t_0 \ioe t \ioe t_0+\fhi(t_0)\big)\\
1 & \big( t_0+\fhi(t_0) \ioe t \ioe T\big).
\end{dcases}
\end{equation}

Les trois branches de $f_2$ et de $f_2'$ se raccordent continûment aux points $t_0$ et $t_0+\fhi(t_0)$ ; pour la fonction $f_2$ au point $t_0$, on utilise l'équation
\[
\fhi(t_0)^2+2t_0\fhi(t_0)-t_0^2-4=0.
\]
La dérivée seconde $f_2''$ ne prend que les valeurs $0$ et $\pm 1$. Comme $\fhi(t_0)\soe \sqrt{2}\soe t_0$, on vérifie que~$f_2(t)$ croît strictement de $f_2(0)=-1$ à $f_2\big(t_0+\fhi(t_0)\big)=1$ lorsque $t$ croît de $0$ à $t_0+\fhi(t_0)$. Ainsi,~$f_2 \in \Lcal$ et $f_2'(t_0)=\fhi(t_0)$, d'où le résultat.

\smallskip

$\bullet$ Si $t_0 > \sqrt{2}$, donc $T-t_0 \soe t_0> \sqrt{2}$, le choix $h=h'=\sqrt{2}$, dans la proposition \ref{190603b}, fournit l'inégalité
\[
\sigma(t_0,T) \ioe G(\sqrt{2},\sqrt{2})=\sqrt{2}.
\]
Si l'on pose $f_3(t)=q(t-t_0-\sqrt{2})$, on a $f_3 \in \Lcal$ et $f_3'(t_0)=\sqrt{2}$, d'où l'égalité $\sigma(t_0,T) =\sqrt{2}$.\fin

\smallskip

Notons que, dans la proposition \ref{190605b}, 

$\bullet$ si $T\ioe 2$, seul le premier cas se produit ;

$\bullet$ si $2 <T\ioe 2\sqrt{2}$, seuls les deux premiers cas se produisent ;

$\bullet$ si $T > 2\sqrt{2}$, seuls les deux derniers cas se produisent.

\subsection{Une fonction de comparaison}\label{200117b}

La fonction suivante joue un rôle central dans l'étude du problème de Landau sur $\Real$. Elle est définie comme la fonction $q : \Real \vers \Real$, de période $4\sqrt{2}$, telle que
\begin{align*}
q(t) &=1-t^2/2 \quad ( \abs{t} \ioe \sqrt{2})\\
q(t+2\sqrt{2}) &=-q(t) \quad (\abs{t} \ioe \sqrt{2}).
\end{align*}

On en déduit que cette dernière relation est valable pour tout $t$, et que $q$ est une fonction paire. On a aussi la relation $q(\sqrt{2}+t)=-q(\sqrt{2}-t)$ pour tout $t$.

\begin{figure}
  \centering
    \includegraphics[width=0.7\textwidth]{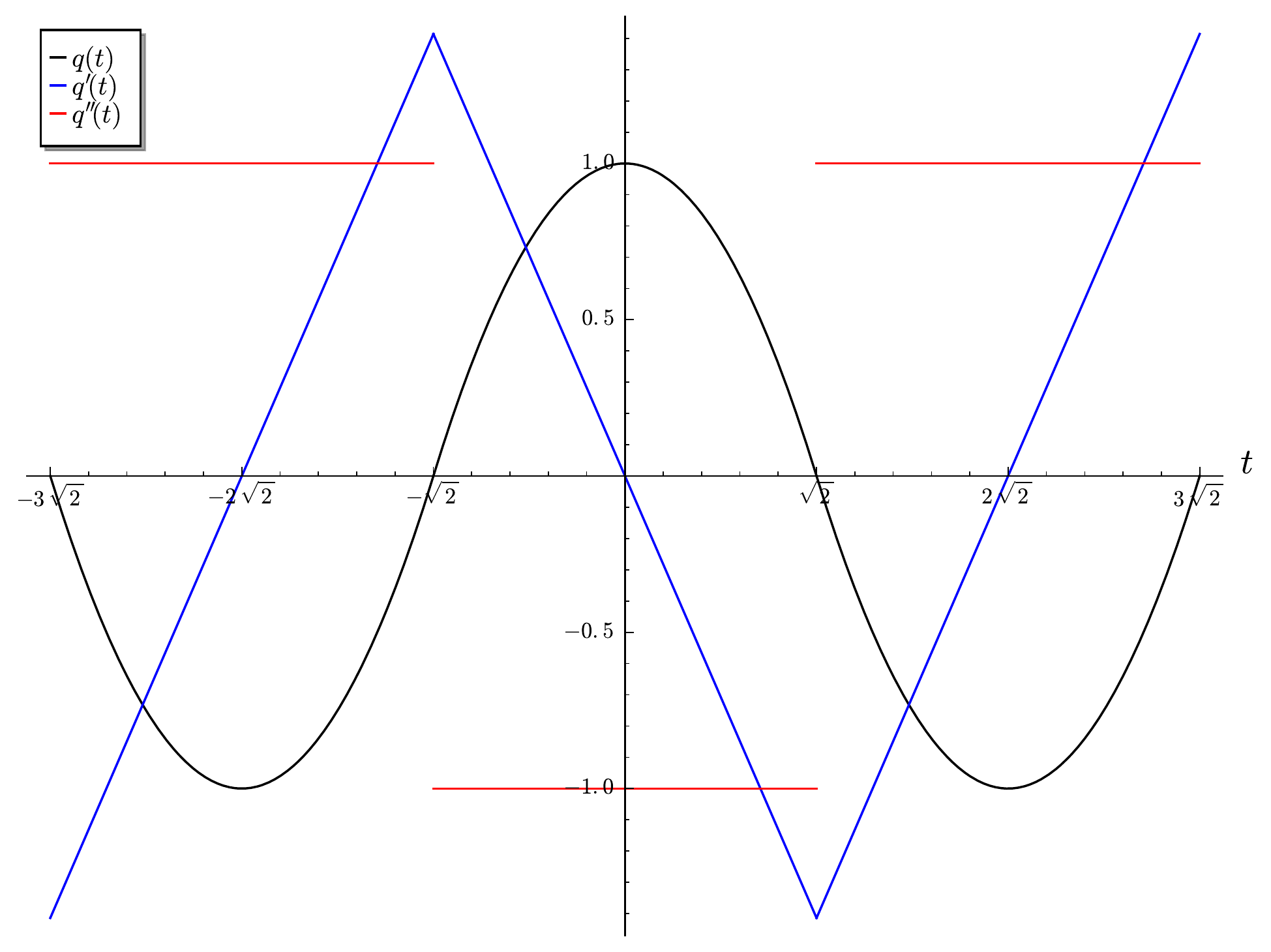}
  \caption{La fonction $q$ et ses dérivées, première et seconde}\label{190410a}
\end{figure}

La figure \ref{190410a} donne une représentation graphique des fonctions $q$, $q'$ et $q''$. La fonction $q$ est dérivable sur $\Real$, et deux fois dérivable sur $\Real \setminus (\sqrt{2}+2\sqrt{2}\,\Int)$, vérifiant
\begin{align*}
q'(t) &=-t \\
q''(t) &= -1,
\end{align*}
pour $\abs{t} \ioe \sqrt{2}$. On a
\[
\norm{q}_{\infty}=\norm{q''}_{\infty}=1 \quad ; \quad \norm{q'}_{\infty}=\sqrt{2},
\]
et la fonction $q$ vérifie l'équation différentielle 
\begin{equation}\label{190625b}
\abs{q}=1-q'^{\,2}/2.
\end{equation}

En revenant au cas général, nous dirons qu'une fonction $\fhi \in \Lcal_n(a,b\, ; \, I)$ est une \emph{fonction de comparaison} pour cet ensemble si
\[
\forall f \in \Lcal_n(a,b\, ; \, I), \; \forall\, t_0,t_1 \in I, \; f(t_0)=\fhi(t_1) \implique \abs{f'(t_0)} \ioe \abs{\fhi'(t_1)}.
\]

Cette notion a été introduite par Kolmogorov (cf. \cite{zbMATH00050290}, p. 281, \cite{zbMATH01369964}, p. 265). Elle est l'objet d'une étude détaillée au chapitre 3 de \cite{zbMATH00047600}.

Nous allons voir que la fonction $q$ est une fonction de comparaison pour $\Lcal_2(1,1\, ; \, \Real)$, mais nous aurons l'usage d'une propriété plus précise.

\begin{prop}\label{200117a}
Soit $T \soe 2\sqrt{2}$, et $t_0$ tel que $\sqrt{2} \, \ioe t_0 \ioe T-\sqrt{2}$. Si $f \in \Lcal_2(T)$ et $t_1 \in \Real$, on a
\[
f(t_0)=q(t_1) \implique \abs{f'(t_0)} \ioe \abs{q'(t_1)}.
\]
\end{prop}
\dem

Supposons d'abord $f(t_0) \soe 0$ et~$f'(t_0)\soe 0$. Soit $t_2$ tel que $-\sqrt{2} \ioe t_2 \ioe 0$, et $ q(t_2)=q(t_1)$. On a donc 
\[
-t_2=q'(t_2)=\abs{q'(t_1)}. 
\]

Si $t_2=0$, on a $f(t_0)=q(t_2)=1$, donc~$f'(t_0)=0=q'(t_2)$. 

Si $-\sqrt{2} \ioe t_2 < 0$, on a
\begin{align*}
1\soe f(t_0-t_2) &=f(t_0)+\int_0^{-t_2}f'(t_0+t)\, dt \\
&=q(t_2)-t_2f'(t_0)+\int_0^{-t_2}\big(f'(t_0+t)-f'(t_0)\big)\, dt\\
&\soe q(t_2) -t_2f'(t_0)+\int_0^{-t_2}(-t)\, dt,
\end{align*}
puisque~$f \in \Lcal_2(T)$. On en déduit
\[
-t_2f'(t_0) \ioe 1-q(t_2)+t_2^2/2=-t_2q'(t_2),
\]
et donc encore $f'(t_0) \ioe q'(t_2)$.

Les cas correspondant aux trois autres possibilités pour $\big(\sgn f(t_0), \sgn f'(t_0)\big)$ se ramènent à celui que nous avons traité, en remplaçant la fonction $f$ par $g$, $t_0$ par $t'_0$, et $t_1$ par $t'_1$, avec,

$\bullet$ pour $(+,-)$, $g(t)=f(T-t)$, $t'_0=T-t_0$ et $t'_1=t_1$ ;

$\bullet$ pour $(-,+)$, $g(t)=-f(T-t)$, $t'_0=T-t_0$ et $t'_1=t_1+2\sqrt{2}$ ;

$\bullet$ pour $(-,-)$, $g(t)=-f(t)$, $t'_0=t_0$ et $t'_1=t_1+2\sqrt{2}$.\fin

\smallskip

La propriété de $q$ d'être une fonction de comparaison pour $\Lcal_2(1,1\, ; \, \Real)$ découle de cette proposition : si $f \in \Lcal_2(1,1\, ; \, \Real)$, si $t_0,\, t_1 \in \Real$ et $f(t_0)=q(t_1)$, la proposition \ref{200117a} appliquée au point $\sqrt{2}$ et à la fonction $g : [0,2\sqrt{2}] \vers \Real$ définie par $g(t)=f(t+t_0-\sqrt{2}\,)$, fournit bien l'inégalité~$\abs{f'(t_0)} \ioe \abs{q'(t_1)}$.

\subsection{Le problème de Landau pour $I=\Real$}\label{200416b}

La solution du problème de Landau sur $\Real$ peut être énoncée sous la forme précisée qui suit.

\begin{prop}\label{190702a}
Pour $f \in \Lcal_2(1,1\, ; \, \Real)$ et $t \in \Real$, on a
\[
\abs{f'(t)} \ioe \sqrt{2\big(1-\abs{f(t)}\big)}.
\]

En particulier,
\begin{equation}\label{190627h}
\sup_{f \in \Lcal_2(1,1\, ;\Real)}\, \norm{f'}_{\infty} =\sqrt{2},
\end{equation}
et cette borne est atteinte pour $f=q$.
\end{prop}
\dem

Comme $\abs{f(t)} \ioe 1$, il existe $t_1 \in \Real$ tel que $f(t)=q(t_1)$. Comme $q$ est une fonction de comparaison pour $\Lcal_2(1,1\, ; \, \Real)$, on a
\begin{equation*}
\abs{f'(t)} \ioe \abs{q'(t_1)}=\sqrt{2\big(1-\abs{q(t_1)}\big)}=\sqrt{2\big(1-\abs{f(t)}\big)},
\end{equation*}
où l'on a utilisé l'équation différentielle \eqref{190625b}.\fin

\smallskip

La fonction $q$ n'est pas l'unique fonction $f \in \Lcal_2(1,1\, ;\Real)$ vérifiant \eqref{190627h}. Un autre exemple (parmi une infinité) est la fonction $t \mapsto q\big(\min(\abs{t},2\sqrt{2}\,)\big)$.

\smallskip

De façon analogue, on déduit de la proposition \ref{200117a} la proposition suivante.

\begin{prop}\label{190703a}
Soit $f \in \Lcal_2(T)$ et $t$ tel que $\sqrt{2} \ioe t \ioe T-\sqrt{2}$. On a
\[
\abs{f'(t)} \ioe \sqrt{2\big(1-\abs{f(t)}\big)}\,.
\]
\end{prop}

\smallskip

En particulier, si $\sqrt{2} \ioe t_0 \ioe T-\sqrt{2}$, on retrouve le fait, vu à la proposition \ref{190605b}, que
\begin{equation*}
\sup_{f \in \Lcal_2(T)}\, \abs{f'(t_0)} =\sqrt{2}.
\end{equation*}

\subsection{Propriétés de prolongement}\label{190624b}

Si $T >0$, $h>0$, et $f \in \Lcal_2(T)$, il est possible que $f$ ne soit pas prolongeable en une fonction appartenant à $\Lcal_2(T+h)$. En revanche, si $0<\eps \ioe 1$ et si $h$ est assez petit, un tel prolongement est possible pour la fonction $(1-\eps)f$, comme le précise la proposition suivante.

\begin{prop}
Soit $T >0$, $h>0$, $0<\eps \ioe 1$, et $f \in \Lcal_2(T)$. Si 
\[
\abs{\thet} \ioe 1 \,\text{ et }\, 0<h \ioe \eps/\big(\abs{\thet}+\sigma_{\infty}(T)\big), 
\]
alors la fonction $g$ définie sur $[0, T+h]$ par 
\begin{equation}\label{190620a}
g(t)=
\begin{dcases}
(1-\eps)f(t) & (0 \ioe t \ioe T)\\
(1-\eps)f(T)+(1-\eps)f'(T)(t-T) +\thet(t-T)^2/2 &(T < t \ioe T+h)
\end{dcases}
\end{equation}
appartient à $\Lcal_2(T+h)$.
\end{prop}
\dem

Les deux branches de la fonction $g$ et de sa dérivée se raccordent au point $T$, et $\abs{g''} \ioe 1$ presque partout. Par conséquent, $g \in \Lcal_2(T+h)$ pourvu que sa valeur absolue soit bornée par $1$ sur~$]T,T+h]$, ce qui sera réalisé si
\[
(1-\eps)\abs{f(T)} +(1-\eps)\abs{f'(T)}h+\abs{\thet}h^2/2 \ioe 1,
\]
ou encore, puisque $\abs{f(T)} \ioe 1$ et $\abs{f'(T)} \ioe \sigma_{\infty}(T)$, si 
\[
0 <h\ioe 2  \,\text{ et }\,  1-\eps+\sigma_{\infty}(T)h+\abs{\thet}h \ioe 1,
\]
\cad si $h\ioe \eps/\big(\abs{\thet}+\sigma_{\infty}(T)\big)$, quantité $\ioe 1/2$.\fin

\smallskip

La proposition suivante décrit un cas où le prolongement de $f$ est possible, non au sens ordinaire, mais, en quelque sorte, {\og à l'intérieur\fg} de sa courbe représentative.

\begin{prop}\label{280608e}
Soit $T >0$ et $f \in \Lcal_2(T)$. On suppose qu'il existe $t_0 \in [0,T[$ tel que~$f'(t_0)=0$ et $f(t_0) <1$. Soit $h$ tel que $0<h\ioe 4\sqrt{1-f(t_0)}$. Alors la fonction $g$ définie sur $[0, T+h]$ par 
\begin{equation}\label{200608a}
g(t)=
\begin{dcases}
f(t) & (0 \ioe t \ioe t_0)\\
f(t_0)+(t-t_0)^2/2&(t_0 < t \ioe t_0+h/4)\\
f(t_0)+h^2/16-(t-t_0-h/2)^2/2&(t_0 +h/4< t \ioe t_0+3h/4)\\
f(t_0)+(t-t_0-h)^2/2&(t_0 +3h/4< t \ioe t_0+h)\\
f(t-h) & (t_0+h<t\ioe T+h)
\end{dcases}
\end{equation}
appartient à $\Lcal_2(T+h)$.
\end{prop}
\dem

Les cinq branches de $g$ et de sa dérivée se raccordent aux points de jonctions des intervalles consécutifs, on a $\abs{g''} \ioe 1$ presque partout, et la {\og bosse\fg}, insérée dans la courbe représentative de $f$, a une hauteur $h^2/16 \ioe 1-f(t_0)$.\fin

\smallskip

Enfin, en adoptant un autre point de vue, on peut se demander à quelle condition un élément~$f \in~\Lcal_2(T)$ est prolongeable en une fonction appartenant à $\Lcal_2(1,1 \, ;\Real)$. La réponse est fournie par la proposition suivante.
\begin{prop}\label{190704a}
Soit $f \in \Lcal_2(T)$. Les quatre assertions suivantes sont équivalentes.

\smallskip

$(i)$ Il existe $g \in \Lcal_2(1,1 \, ;\Real)$ dont la restriction à $[0,T]$ coïncide avec~$f$, et telle que $g'$ est à support compact ;

$(ii)$ Il existe $g \in \Lcal_2(1,1 \, ;\Real)$ dont la restriction à $[0,T]$ coïncide avec~$f$ ;

$(iii)$ Pour tout $t\in [0,T]$, on a $\abs{f'(t)} \ioe \sqrt{2\big(1-\abs{f(t)}\big)}\,$ ;

$(iv)$ On a $\abs{f'(0)} \ioe \sqrt{2\big(1-\abs{f(0)}\big)}\,$ et $\abs{f'(T)} \ioe \sqrt{2\big(1-\abs{f(T)}\big)}\,$.

\end{prop}
\dem

Les implications $(i) \implique (ii)$ et $(iii) \implique (iv)$ sont triviales, et l'implication $(ii) \implique (iii)$ découle de la proposition \ref{190702a}. Il reste à démontrer que $(iv)$ entraîne $(i)$.

\smallskip

Si la condition $(iv)$ est remplie, définissons une fonction $g: \Real \vers \Real$ par
\begin{equation}\label{190624a}
g(t)=
\begin{dcases}
f(0)-\sgn f'(0)\cdot f'(0)^2/2& (t\ioe-\abs{f'(0)})\\
f(0) +f'(0)t+\sgn f'(0)\cdot t^2/2&(-\abs{f'(0)} < t \ioe 0)\\
f(t) & (0 < t \ioe T)\\
f(T)+f'(T)(t-T)-\sgn f'(T)\cdot (t-T)^2/2&(T < t \ioe T+\abs{f'(T)})\\
f(T)+\sgn f'(T)\cdot f'(T)^2/2& (t>T+\abs{f'(T)})
\end{dcases}
\end{equation}

Les branches de $g$ se raccordent, et de même pour $g'$, aux points de jonctions. La fonction $g$ est continûment dérivable, et $g'$ est lipschitzienne de rapport $1$, et à support compact. De plus, la deuxième et la quatrième branche sont monotones, les sens de variation étant donnés par les signes de $f'(0)$ et $f'(T)$, respectivement. Pour vérifier que~$g\in~\Lcal_2(1,1 \, ;\Real)$, il reste donc seulement à vérifier que les constantes définissant la première et la cinquième branche de $g$ ont des valeurs absolues $\ioe 1$, ce qui résulte de $(iv)$.\fin 

\smallskip

Observons que l'équivalence entre $(iii)$ et $(iv)$ est une sorte de {\og principe du maximum\fg} dans l'ensemble $\Lcal_2(T)$. 

La proposition suivante est issue du Lemma 1, p. 60 de \cite{zbMATH01891483}

\begin{prop}\label{190627b}
Soit $T > 2\sqrt{2}$ et $f \in \Lcal_2(T)$. Alors il existe $g \in \Lcal_2(1,1 \, ;\Real)$, dont la dérivée est à support compact, et telle que $f(t)=g(t)$ pour tout $t \in [\sqrt{2},T-\sqrt{2}\,]$. 
\end{prop}
\dem

D'après la proposition \ref{190703a}, la fonction $f$ vérifie 
\[
\abs{f'(t)} \ioe \sqrt{2\big(1-\abs{f(t)}\big)} \quad (\sqrt{2}\ioe t \ioe T-\sqrt{2}\,).
\]

Le résultat découle donc de la proposition \ref{190704a}, avec le segment $[0,T]$ remplacé par le segment~$[\sqrt{2},T-\sqrt{2}\,]$.\fin

\subsection{Le problème extrémal pour la vitesse moyenne}\label{190406b}

Si $f \in \Lcal_2(T)$, la quantité
\[
\Phi_T(f)=\int_0^T \abs{f'(t)}\, dt
\]
est la variation totale de $f$. Dans l'interprétation cinématique, c'est la distance totale parcourue par la particule lors de ses allées et venues dans le segment $[-1,1]$, pendant l'intervalle de temps~$0 \ioe t \ioe T$. Quant à $\Phi_T(f)/T$, c'est la valeur moyenne de la valeur absolue de la vitesse de la particule.

La proposition \ref{190615a} s'applique à $\Phi_T$, qui est convexe et continue pour la topologie de $(W,N_1)$. Pour $T>0$, nous noterons
\[
\sigma_1(T)=\sup_{f \in \Lcal_2(T)} \Phi_T(f).
\]
Par convention, $\sigma_1(0)=0$.

L'étude de $\sigma_1(T)$ et des fonctions extrémales correspondantes a été notamment abordée, sous une forme plus générale, par Bojanov et Naidenov (cf. \cite{zbMATH01369964} et \cite{zbMATH01891483}), puis Naidenov (cf. \cite{zbMATH01970534}). Le présent paragraphe est une introduction à ces travaux.

\subsubsection{Observations préalables}

Nous utiliserons les trois faits suivants.

$\bullet$ Si $f \in \Lcal_2(T)$ et si $f'$ ne s'annule pas dans l'intervalle $]\alpha,\beta[ \,\subset [0,T]$, alors
\begin{equation*}
\int_\alpha^\beta \abs{f'(t)} \, dt =\abs{\int_{\alpha}^{\beta} f'(t) \, dt }=\abs{f(\beta)-f(\alpha)}\ioe 2.
\end{equation*}

$\bullet$ Si $f \in \Lcal_2(T)$ et si $f'(t_0)=0$, où $0\ioe t_0 \ioe T$, alors
\begin{equation*}
\abs{f'(t) }\ioe \abs{t-t_0} \quad (0 \ioe t \ioe T).
\end{equation*}

$\bullet$ Si $\alpha \ioe t_0\ioe \beta$, l'intégrale 
\begin{equation*}
\int_{\alpha}^{\beta} \abs{t-t_0} \, dt =\frac{(t_0-\alpha)^2+(\beta-t_0)^2}2
\end{equation*}
est une fonction de la variable $t_0$, décroissante sur $[\alpha,(\alpha+\beta)/2]$, croissante sur $[(\alpha+\beta)/2,\beta]$. Son maximum est atteint aux points $\alpha$ et $\beta$, et vaut $(\beta-\alpha)^2/2$.

\subsubsection{Valeur de $\sigma_1(T)$ pour $T \ioe 4$}

\begin{prop}\label{190626b}
Pour $0 < T \ioe 2$, on a $\sigma_1(T)=2$.
\end{prop}
\dem

Supposons $0<T \ioe 2$. La fonction 
\[
\tau(t)=-1+(2/T-T/2)t+t^2/2 \quad (0 \ioe t \ioe T),
\]
vérifie $\tau''=1$, $\tau(0)=-1$, $\tau(T)=1$, et
\[
\tau'(t)=t+2/T-T/2 \soe t \soe 0 \quad (0 \ioe t \ioe T).
\]

Elle appartient donc à $\Lcal_2(T)$ (on a d'ailleurs $\tau(t)=-f(T-t)$, si $f$ est définie par \eqref{190620b}). Comme
\[
\int_0^T \abs{\tau'(t)} \, dt = \int_0^T \tau'(t) \, dt =\tau(T)-\tau(0)=2,
\]
on a $\sigma_1(T) \soe 2$.

\smallskip

Maintenant, soit $f \in \Lcal_2(T)$. Examinons deux cas.

\underline{Premier cas :} il existe $t_0 \in [0,T]$ tel que $f'(t_0)=0$.

On a alors
\[
\int_0^T \abs{f'(t)} \, dt \ioe \int_0^{T} \abs{t-t_0} \, dt \ioe T^2/2\ioe 2.
\]

\underline{Second cas :} la dérivée $f'$ ne s'annule pas sur $[0,T]$. On a donc $\Phi_T(f) \ioe 2$.

Par conséquent, $\sigma_1(T) \ioe 2$.\fin

\begin{prop}\label{200608b}
Pour $2 \ioe T\ioe 4$, on a
\[
\sigma_1(T)=T^2/2-2T+4.
\]
\end{prop}
\dem

Nous explicitons dans ce cas particulier la remarque de Bojanov et Naidenov, p. 68 de \cite{zbMATH01891483}.

Posons
\begin{equation}\label{200608c}
\fhi(t)=-1+(t-2)^2/2 \quad (0 \ioe t \ioe T).
\end{equation}

On a $\fhi \in \Lcal_2(T)$, et
\[
\Phi_T(\fhi)=\int_0^T \abs{t-2} \, dt=2+(T-2)^2/2=T^2/2-2T+4, 
\]
donc $\sigma_1(T)\soe T^2/2-2T+4$.

\smallskip

Pour démontrer l'inégalité en sens inverse, considérons $f \in \Lcal_2(T)$. Examinons deux cas.

\smallskip

\underline{Premier cas} : la dérivée $f'$ ne s'annule pas sur $[0,T]$. On a alors $\Phi_T(f) \ioe 2 < 2+(T-2)^2/2$.

\underline{Second cas} : la dérivée $f'$ s'annule en un point $t_0 \in [0,T]$ ; on a donc $\abs{f'(t)} \ioe \abs{t-t_0}$ pour tout $t \in [0,T]$.

\underline{Premier sous-cas} : $t_0 \soe 2$. On a alors
\begin{align*}
\Phi_T(f) &=\int_0^2 \abs{f'(t)} \, dt +\int_{2}^T \abs{f'(t)} \, dt\\
&\ioe \sigma_1(2)+\int_{2}^T \abs{t-t_0} \, dt\\
&\ioe 2+ (T-2)^2/2.
\end{align*}

\underline{Deuxième sous-cas} : $t_0 \ioe T-2$. Ce sous-cas se ramène au précédent en remplaçant la fonction~$f$ par la fonction $t \mapsto f(T-t)$.

\underline{Troisième sous-cas} : $T-2 < t_0 < 2$. On a alors
\begin{equation*}
\Phi_T(f) \ioe  \int_{0}^{T} \abs{t-t_0} \, dt\ioe \int_{0}^{T} \abs{t-2} \, dt=2+ (T-2)^2/2.\fine
\end{equation*}

\subsubsection{Sous-additivité}

\begin{prop}\label{200427a}
La fonction $\sigma_1$ est sous-additive : pour~$T_1,T_2 > 0$, on a
\[
\sigma_1(T_1+T_2) \ioe \sigma_1(T_1)+\sigma_1(T_2).
\]
\end{prop}
\dem

Pour $T=T_1+T_2$ et $f \in \Lcal_2(T)$, on a
\begin{equation*}
\Phi_T(f)= \int_0^{T _1}\abs{f'(t)}\, dt +\int_0^{T_2} \abs{f'(T_1+t)}\, dt\ioe \sigma_1(T_1)+\sigma_1(T_2).\fine
\end{equation*}

\smallskip

En écrivant~$T=2\floor{T/2} +2\set{T/2}$, on en déduit une première majoration de $\sigma_1(T)$ :
\begin{equation}\label{190616a}
\sigma_1(T)  \ioe \floor{T/2} \sigma_1(2) +\sigma_1(2\set{T/2})
\ioe 2\floor{T/2} +2
\ioe T+2 \quad(T>0).
\end{equation}

\smallskip

Pour $4\ioe T \ioe 8$, on déduit des propositions \ref{200608b} et \ref{200427a} la majoration
\begin{equation}\label{200629a}
\sigma_1(T)  \ioe 2\sigma_1(T/2)=4+(T-4)^2/4 \quad(4\ioe T \ioe 8).
\end{equation}

\begin{prop}\label{190627a}
La limite
\[
\lim_{T \vers \infty} \sigma_1(T)/T
\]
existe, et est égale à
\[
\inf_{T >0} \sigma_1(T)/T.
\]
\end{prop}
\dem

D'après \eqref{190616a}, la fonction $\sigma_1$ est bornée sur tout intervalle $]0,A]$ (par $A +2$). Le résultat découle alors de la sous-additivité de $\sigma_1$ et de la version continue du lemme de Fekete donnée par Hille et Phillips dans \cite{zbMATH03127846}, Theorem 7.6.1, p. 244.\fin

\smallskip

La proposition \ref{190615b} ci-dessous prouvera que cette limite vaut $1/\sqrt{2}$.

\subsubsection{Continuité}

\begin{prop}\label{200608c}
La fonction $\sigma_1$ est continue sur $]0,\infty[$.
\end{prop}
\dem

Nous allons montrer les inégalités
\begin{align}
\sigma_1(T+h) &\ioe \sigma_1(T) + h\sigma_{\infty}(T) \quad ( T>0, \, h>0)\label{190617c}\\
\sigma_1(T+h) &\soe \sigma_1(T) - h\sigma_{\infty}(T)\sigma_1(T) \quad ( T>0, \, 0<h \ioe 1/\sigma_{\infty}(T))\label{190617d}
\end{align}

Compte tenu de \eqref{190616a}, comme $\sigma_{\infty}(T)$ est une fonction décroissante de $T$, il en résultera que, pour tous $\alpha$, $A$ tels que $0 < \alpha <A$, nous aurons
\[
\abs{\sigma_1(T+h) - \sigma_1(T)} \ioe Ch \quad (\alpha \ioe T \ioe A,\; 0 <h \ioe 1/\sigma_{\infty}(\alpha)\,),
\]
avec $C=(A+2)\sigma_{\infty}(\alpha)$, d'où la continuité annoncée ; $\sigma_1$ est même lipschitzienne sur tout segment de $]0,\infty[$.

\smallskip

Démontrons \eqref{190617c}. Pour~$T,h> 0$, et $f \in \Lcal_2(T+h)$, on a

\begin{equation*}
\int_{T}^{T+h} \abs{f'(t)}\, dt\ioe h\,\sigma_{\infty}(T+h)\ioe h\,\sigma_{\infty}(T),
\end{equation*}
donc 
\begin{equation*}
\Phi_{T+h}(f)= \int_0^{T }\abs{f'(t)}\, dt +\int_T^{T+h} \abs{f'(t)}\, dt\ioe \sigma_1(T)+h\sigma_{\infty}(T).
\end{equation*}

\smallskip

Démontrons \eqref{190617d}. Soit $T> 0$ et $f \in \Lcal_2(T)$ telle que $\Phi_{T}(f)=\sigma_1(T)$. Si $0<h \ioe 1/\sigma_{\infty}(T)$, posons $\eps=h\sigma_{\infty}(T)$, et considérons la fonction $g$ définie par~\eqref{190620a}, avec $\thet=0$. Elle appartient à~$\Lcal_2(T+h)$, et
\begin{align*}
\sigma_1(T+h) &\soe\int_0^{T} \abs{g'(t)} \, dt = (1-\eps)\int_0^{T }\abs{f'(t)}\, dt \\
&=(1-\eps)\sigma_1(T)=\sigma_1(T)-h\sigma_{\infty}(T)\sigma_1(T).\fine
\end{align*}

\subsubsection{Monotonie}

\begin{prop}
La fonction $\sigma_1$ est strictement croissante sur $[2,\infty[$.
\end{prop}
\dem

La proposition \ref{200608b} permet de se ramener à l'intervalle $[4,\infty[$. On a d'abord la minoration~$\sigma_1(T) \soe 4$ pour tout $T \soe 4$, car la fonction $f_0$ définie par 
\[
f_0(t)=
\begin{dcases}
-1+(t-2)^2/2 &(0\ioe t \ioe 2)\\
0 &(2 < t <T-2)\\
-1+(t-T+2)^2/2 & (T-2\ioe t \ioe T)
\end{dcases}
\]
appartient à $\Lcal_2(T)$ et vérifie $\Phi_{T}(f_0)=4$. Ensuite, soit $f \in \Lcal_2(T)$ telle que $\Phi_{T}(f)=\sigma_1(T)$. Comme $\Phi_{T}(f) >2$, il existe $t_0 \in \, ]0,T[$ tel que $f'(t_0)=0$. Quitte à remplacer $f$ par $-f$, on peut supposer que $f(t_0) <1$. Pour $0<h\ioe 4\sqrt{1-f(t_0)}$, la fonction $g$ définie dans la proposition \ref{280608e} appartient à $\Lcal_2(T+h)$ et on a
\[
\sigma_1(T+h)\soe \Phi_{T+h}(g)=\Phi_{T}(f)+h^2/8 > \sigma_1(T)
\]
La croissance stricte de $\sigma_1$ sur $[4,\infty[$ résulte de cette croissance stricte locale, et de la continuité de $\sigma_1$ (proposition~\ref{200608c}).\fin

\subsubsection{Comportement asymptotique de $\sigma_1(T)$ lorsque $T$ tend vers l'infini}

Bojanov et Naidenov ont remarqué  dans \cite{zbMATH01891483} (Theorem 2, p. 68), que l'on peut déterminer la valeur exacte de $\sigma_1(T)$ pour une suite de valeurs de $T$ tendant vers l'infini. La clé de ce calcul est la proposition suivante.

\begin{prop}\label{190626a}
Soit $f \in \Lcal_2(1,1 \, ;\Real)$ telle que $f'$ est à support compact. On a
\[
\int_{-\sqrt{2}}^{\sqrt{2}}\abs{f'(t)} \, dt \ioe 2.
\]
\end{prop}
\dem

Nous suivons ici la démonstration du Theorem 1, p. 268 de \cite{zbMATH01369964}.

Posons
\[
\fhi(u)=\int _{-\sqrt{2}}^{\sqrt{2}}\abs{f'(u+t)} \, dt \quad (u \in \Real).
\]
Il s'agit de montrer que $\fhi(0)\ioe 2$.

La fonction $\fhi$ est continûment dérivable sur $\Real$, et à support compact. Elle atteint donc son maximum global, en un point $u_0\in \Real$, où sa dérivée s'annule :
\[
0=\fhi'(u_0)=\abs{f'(u_0+\sqrt{2})}-\abs{f'(u_0-\sqrt{2})}.
\]

Posons $g(t)=f(u_0+t)$ pour $t \in \Real$. La fonction $g$ appartient à $\Lcal_2(1,1 \, ;\Real)$ et 
\begin{equation}\label{190626c}
\abs{g'(\sqrt{2})}=\abs{g'(-\sqrt{2})}.
\end{equation}

Comme pour la proposition \ref{190626b}, examinons deux cas.

\smallskip

\underline{Premier cas :} il existe $t_0 \in [-\sqrt{2},\sqrt{2}]$ tel que $g'(t_0)=0$. 

Définissons alors une fonction $Q$ sur $[-\sqrt{2},\sqrt{2}]$ par
\[
Q(t)=
\begin{dcases}
\lvert g'(t_0-t-\sqrt{2})\rvert  &(-\sqrt{2} \ioe  t \ioe t_0)\\
\lvert g'(t_0-t+\sqrt{2})\rvert  &\quad (t_0 \ioe t \ioe \sqrt{2}).
\end{dcases}
\]

Les deux branches se raccordent au point $t_0$, d'après \eqref{190626c}. La fonction $Q$ est donc, comme $g'$, lipschitzienne de rapport $1$. De plus $Q(\pm\sqrt{2})=\abs{g'(t_0)}=0$. On a donc
\[
Q(t) \ioe \sqrt{2}-\abs{t} \quad ( -\sqrt{2} \ioe t \ioe \sqrt{2}).
\]

Par conséquent,
\begin{align*}
\fhi(0) & \ioe \fhi(u_0) =\int_{-\sqrt{2}}^{\sqrt{2}}\abs{g'(t)} \, dt\\
&=\int_{-\sqrt{2}}^{t_0}\abs{g'(t)} \, dt +\int_{t_0}^{\sqrt{2}}\abs{g'(t)} \, dt= \int_{-\sqrt{2}}^{\sqrt{2}} Q(t) \, dt\\
&\ioe \int_{-\sqrt{2}}^{\sqrt{2}} \big(\sqrt{2}-\abs{t}\big)  \, dt =2.
\end{align*}

\smallskip

\underline{Second cas :} la dérivée $g'$ ne s'annule pas sur $[-\sqrt{2},\sqrt{2}]$.

Elle est donc de signe constant, et
\begin{equation*}
\fhi(0) \ioe \int_{-\sqrt{2}}^{\sqrt{2}}\abs{g'(t)} \, dt   =\abs{\int_{-\sqrt{2}}^{\sqrt{2}}g'(t) \, dt }=\abs{g(\sqrt{2})-g(-\sqrt{2})}\ioe 2.\fine
\end{equation*}

La borne $2$ de la proposition \ref{190626a} est optimale ; elle est atteinte, par exemple, pour 
\[
f(t)=
\begin{dcases}
q(t-\sqrt{2}\,) & ( \abs{t} \ioe \sqrt{2})\\
\sgn(t) & ( \abs{t} \soe \sqrt{2}).
\end{dcases}
\]

\begin{prop}\label{190615b}
Pour $N\in \Nat$, on a 
\[
 \sigma_1\big(2N\sqrt{2}+4) =2N+4.
\]
\end{prop}
\dem

Le résultat est vrai pour $N=0$, d'après la proposition \ref{200608b}.

Supposons $N\soe 1$. Soit $T=2N\sqrt{2}+4$ et $f \in \Lcal_2(T)$. D'après la proposition \ref{190627b}, $f$ coïncide sur le segment~$[\sqrt{2},T-\sqrt{2}]$ avec une fonction $g\in \Lcal_2(1,1 \, ;\Real)$ telle que $g'$ est à support compact. 

On a
\begin{align*}
\int_{0}^{T}\abs{f'(t)} \, dt &=\int_{0}^{2}\abs{f'(t)} \, dt +\int_{2}^{2N\sqrt{2}+2}\abs{g'(t)} \, dt +\int_{T-2}^{T}\abs{f'(t)} \, dt \\
&=\int_{0}^{2}\abs{f'(t)} \, dt +\int_{T-2}^{T}\abs{f'(t)} \, dt + \sum_{1\ioe k\ioe N}\int_{-\sqrt{2}}^{\sqrt{2}}\abs{g'\big((2k-1)\sqrt{2}+t+2\big)} \, dt \\
&\ioe 2\sigma_1(2)+2N = 2N+4.
\end{align*}
d'après la proposition \ref{190626a}.

Cette majoration est atteinte lorsque 
\begin{equation}\label{200613b}
f(t)=\begin{dcases}
1-(t-2)^2/2& (0 \ioe t \ioe 2)\\
q(t-2)& (2 \ioe t \ioe 2N\sqrt{2}+2)\\
(-1)^N\big(1-(t-2N\sqrt{2}-2)^2/2\big)& (2N\sqrt{2}+2 \ioe  2N\sqrt{2}+4).
\end{dcases}
\end{equation}

\begin{prop}\label{190627e}
Pour tout nombre réel positif $T$, on a l'encadrement
\[
T/\sqrt{2} \ioe \sigma_1(T) \ioe T/\sqrt{2} \,+5.
\]
\end{prop}
\dem

La minoration résulte de la proposition \ref{190627a} et de la valeur $1/\sqrt{2}$ de la limite de $\sigma_1(U)/U$ quand $U$ tend vers l'infini, valeur déterminée grâce à la proposition \ref{190615b}.

Pour la majoration, d'une part, si $T <4\sqrt{2}$, on a $\sigma_1(T) \ioe 5$, d'après \eqref{200629a}. D'autre part, si~$T \soe 4\sqrt{2}$, on écrit~$T=2N\sqrt{2} +y$, où $N$ est un nombre entier $\soe 2$ et $0 \ioe y  < 2\sqrt{2}$, d'où
\begin{align*}
\sigma_1(T)&=\sigma_1\big(2(N-2)\sqrt{2} +4+z\big)\expli{avec $z=y+4\sqrt{2}-4 \in[4\sqrt{2}-4,6\sqrt{2}-4[\,$}\\
& \ioe \sigma_1\big(2(N-2)\sqrt{2} +4\big) +\sigma_1(z)\expli{d'après la proposition \ref{200427a}}\\
& \ioe 2(N-2) +4 +5\expli{d'après la proposition \ref{190615b}, et \eqref{200629a} donnant $\sigma_1(z) \ioe 5$}\\
& = 2N+5 \ioe T/\sqrt{2} \,+5.\fine
\end{align*}

\smallskip

Il serait intéressant de déterminer $\sigma_1(T)$ pour tout $T$.

\subsection{Formulaire}\label{190528b}

Revenons maintenant à l'ensemble $\Lcal_2(a,b\, ; \, I)$, où $a$ et $b$ sont des nombres réels positifs quelconques, et $I$ un segment de longueur $T$. Posons
\begin{align*}
\sigma(a,b\, ; \, t_0,T) &= \sup_{f \in \Lcal_2(a,b\, ; \, I)} \abs{f'(t_0)} \\
\sigma_{\infty}(a,b\, ; \, T) &= \sup_{f \in \Lcal_2(a,b\, ; \, I)} \norm{f'}_{\infty} \\
\sigma_1(a,b\, ; \, T) &= \sup_{f \in \Lcal_2(a,b\, ; \, I)} \int_{0}^{T}\abs{f'(t)} \, dt .
\end{align*}

Comme nous l'avons vu au \S\ref{190322a}, en posant $f_1(t)=(1/a)f(t\sqrt{a/b}\,)$, la transformation $f \mapsto f_1$ applique $\Lcal_2(a,b\, ; \, [0,T])$ sur $\Lcal_2(1,1\, ; \, [0,T\sqrt{b/a}\,])$. On a $f(t)=af_1(t\sqrt{b/a}\,)$, donc
\begin{align*}
f'(t_0) &= \sqrt{ab}\,f'_1(t_0\sqrt{b/a}\,) \quad (0 \ioe t_0 \ioe T)\\
\norm{f'}_{\infty}  &= \sqrt{ab}\, \norm{f'_1}_{\infty}\\
\int_{0}^{T}\abs{f'(t)} \, dt  &=a\int_{0}^{T\sqrt{b/a}}\abs{f_1'(t)} \, dt.
\end{align*}

\smallskip

On déduit de la proposition \ref{190605d} les formules
 \begin{equation}\label{200330a}
\sigma_{\infty}(a,b\, ; \, T)=\sqrt{ab}\cdot\sigma_{\infty}(1,1\, ; \, T\sqrt{b/a}\,)=
\begin{dcases}
2a/T+bT/2 &(T \ioe 2\sqrt{a/b}\,) \\
2 \sqrt{ab}& (T\soe 2\sqrt{a/b}\,).
\end{dcases}
\end{equation}

Comme l'a observé Hadamard (cf. \cite{h14}, p. 70), dans le premier de ces deux cas, on a la majoration~$\sigma_{\infty}(a,b\, ; \, T)\ioe 4a/T$.

\smallskip

La proposition \ref{190627g} fournit la relation
\[
\sup_{f \in \Lcal_2(a,b \, ; \,[0,\infty[)}\, \norm{f'}_{\infty} =2\sqrt{ab}.
\]

\smallskip

Sous la condition $0\ioe t_0\ioe T/2$, on déduit de la proposition \ref{190605b} les formules
\begin{align*}
\sigma(a,b\, ; \, t_0,T) &=\sqrt{ab}\cdot\sigma_{\infty}(1,1\, ; \, t_0\sqrt{b/a}\,\, ; \, T\sqrt{b/a}\,)\\
&=\begin{dcases}
2a/T+bT/2-bt_0(T-t_0)/T &(t_0 \ioe \sqrt{2a/b} \,\text{ et }\, T \ioe \sqrt{2t_0^2+4a/b}\,)\\
\sqrt{2t_0^2b^2+4ab}\,-bt_0 &(t_0 \ioe \sqrt{2a/b} \,\text{ et }\, T > \sqrt{2t_0^2+4a/b}\,)\\
\sqrt{2ab}&(t_0 > \sqrt{2a/b}).
\end{dcases}
\end{align*}

\smallskip

L'égalité \eqref{190627h} fournit l'égalité
\[
\sup_{f \in \Lcal_2(a,b \, ; \, \Real)}\, \norm{f'}_{\infty} =\sqrt{2ab}.
\]

\smallskip

Enfin, on déduit des propositions \ref{190626b}, \ref{200608b} et \ref{190615b} les formules
\begin{align*}
\sigma_1(a,b\, ; \, T) &=a \cdot\sigma_{1}(1,1\, ; \, T\sqrt{b/a}\,)\\
&=
\begin{dcases}
2a &(T \ioe 2\sqrt{a/b}\,) \\
bT^2/2-2T\sqrt{ab} \, +4a &(2\sqrt{a/b}\, <T \ioe 4\sqrt{a/b}\,)\\
(2N+4)a& (T= 2N\sqrt{2a/b}\,+4\sqrt{a/b},\, N \in \Nat),
\end{dcases}
\end{align*}
et de la proposition \ref{190627e} l'encadrement
\[
T\sqrt{ab/2}\, \ioe \sigma_1(a,b\, ; \, T) \ioe T\sqrt{ab/2}\,  +5a \quad (T>0).
\]

\section{Introduction au cas général}\label{200311a}

\subsection{Splines d'Euler}\label{200416c}

Chilov a conjecturé, et Kolmogorov a démontré en 1939, qu'une certaine suite de fonctions, analogue à la suite des fonctions de Bernoulli, fournit la solution du problème de maximiser les normes des dérivées intermédiaires $\|f^{(k)}\|_{\infty}$ $(0 <k<n)$ des fonctions de $\Lcal_n(a,b\, ; \Real)$. Voici cette construction, essentiellement sous la forme que lui a donnée Schoenberg dans son article \cite{zbMATH03411359}.   

\smallskip

On pose 
\begin{align*}
e_0(x) &=
\begin{dcases}
(-1)^{\lfloor x\rfloor} & (x \notin \Int)\\
0 & (x \in \Int)
\end{dcases}\\
e_n(x) &=\int_0^xne_{n-1}(t) \, dt - c_n\quad ( n \in \Nat^*, \, x \in \Real),
\end{align*}
où la constante $c_n$ est choisie pour que l'intégrale $\int_0^2 e_n(t) \, dt$ soit nulle :
\[
c_n=\demi\int_0^2 \Big(\int_0^xne_{n-1}(t) \, dt\Big) \, dx=\int_0^2 ne_{n-1}(t)(1-t/2) \, dt.
\]

La fonction $e_0$ est normalisée ($e_0(x)=\big(e_0(x+0)+e_0(x-0)\big)/2$), et impaire.

\smallskip

On montre alors, par récurrence sur $n \in \Nat$, que

\smallskip

$\bullet$ $e_n$ est de classe $\Ccal^{n-1}$ pour $n \soe 1$ ;

$\bullet$ $e_n$ est de période $2$ ;

$\bullet$ $e_n$ coïncide sur $]0,1[$ avec un polynôme de degré $n$, disons $E_n$ ;

$\bullet$ $c_n=\demi\int_0^1 ne_{n-1}(t) \, dt$ pour $n \soe 1$ ;

$\bullet$ $e_n(x+1)=-e_n(x)$ pour $x \in \Real$ ;

$\bullet$ $e_n(1-x)=(-1)^ne_n(x)$ pour $x \in \Real$ ;

$\bullet$ $c_n=0$ pour $n$ pair, $n \soe 2$ ;

$\bullet$ $e_n$ est impaire si $n$ est pair, et paire si $n$ est impair.

\smallskip

Notons la relation
\begin{equation}\label{200317b}
e_n^{(k)}=\frac{n!}{(n-k)!}\, e_{n-k} \quad (0 \ioe k \ioe n)
\end{equation}
(pour $k=n$, cette relation est valable sur $\Real\setminus \Int$).

\smallskip

Les premiers polynômes $E_n$ sont
\begin{align*}
E_0(x) &= 1 \\
E_1(x) &= x- 1/2 \\
E_2(x) &= x^2-x \\
E_3(x) &= x^3-3 x^2/2 +  1/4 \\
E_4(x) &=  x^4-2x^3+x.
\end{align*}

Cette suite vérifie, comme les fonctions $e_n$, la relation $E'_n=nE_{n-1}$. C'est donc une suite de polynômes d'Appell (cf. \cite{zbMATH02708259}). Par conséquent, on a l'identité entre séries formelles en deux variables (éléments de $\Rat[x][[z]]$),
\[
\sum_{n\soe 0} E_n(x)\frac{z^n}{n!}=E(z)e^{zx},
\]
où $E(z)$ est une série formelle à déterminer, en la variable $z$. Pour cela, on remarque que la relation $e_n(x)+e_n(x+1)=0$ entraîne
\[
E_n(0)+E_n(1)= 2[n=0] \quad ( n \in \Nat),
\]
donc
\[
E(z)+E(z)e^z=2,\]
\cad
\begin{equation}\label{200307a}
\sum_{n\soe 0} E_n(x)\frac{z^n}{n!}=\frac{2e^{zx}}{e^z+1}\cdotp
\end{equation}

Les $E_n$ sont les \emph{polynômes d'Euler} (cf. \cite{zbMATH02596069}, pp. 23-29). 

L'identité
\[
\frac{2e^{zx}}{e^z+1}=\frac{e^{z(x-1/2)}}{\cosh z/2}\cdotp
\]
entraîne
\begin{equation}\label{200307b}
E_n(x)=\sum_{0\ioe \nu \ioe n}\binom{n}{\nu} \frac{E_{\nu}}{2^{\nu}}(x-1/2)^{n-\nu},
\end{equation}
où les $E_{\nu}$ sont les \emph{nombres d'Euler}\footnote{Comme pour les nombres et polynômes de Bernoulli, on adopte la même notation pour les nombres et polynômes d'Euler, ce qui ne doit mener à aucune confusion.} (suite A122045 de \cite{oeis}), définis par leur série génératrice exponentielle
\begin{equation*}
\sum_{\nu \soe 0} E_{\nu}\frac{z^{\nu}}{\nu!}=\frac{1}{\cosh z}\cdotp
\end{equation*}

Les nombres d'Euler sont entiers. Ceux d'indice impair sont nuls, et ceux d'indice pair ont des signes alternés. Les premiers nombres d'Euler non nuls sont
\[
E_0=1 \quad ; \quad E_2=-1 \quad ; \quad E_4=5 \quad ; \quad E_6=-61 \quad ; \quad E_8=1385.
\]

\smallskip

Rappelons qu'une fonction $f$, à valeurs réelles, définie sur un segment $[\alpha,\beta]$, est une \emph{spline} de degré $k \in \Nat^*$, si elle est de classe $\Ccal^{k-1}$, et s'il existe une subdivision
\[
t_0=\alpha < t_1 < \dots < t_j=\beta
\]
et $j$ polynômes de degrés $\ioe k$, disons $P_1$,..., $P_j$, tels que
\begin{equation}\label{190603a}
f(t)=P_i(t) \quad (t_{i-1} <t < t_i, \,1\ioe i \ioe j).
\end{equation}
On suppose en général que cette décomposition est minimale, \cad que $P_i$ et $P_{i+1}$ sont distincts pour tout $i$. Par convention, une spline de degré $0$ est simplement une fonction en escalier sur $[\alpha,\beta]$. On peut étendre la notion de spline de degré $k$ au cas d'une fonction définie sur toute la droite réelle : sa restriction à tout segment doit être une spline de degré $k$ sur ce segment.

Il résulte de ce qui précède que, pour tout $n \in \Nat$, la fonction $e_n$ est une spline de degré $n$ sur~$\Real$. \smallskip

Pour $n \in \Nat$, posons
\[
r_n=\sup_{x \in \Real} \abs{e_n(x)}.
\]

Compte tenu des relations $e_n(x+1)=-e_n(x)$ et $e_n(1-x)=(-1)^ne_n(x)$, on a
\[
r_n=\sup_{0 \ioe x \ioe 1/2} \abs{e_n(x)}.
\]

On montre par récurrence sur $n \soe 1$ que 

$\bullet$ $\abs{e_n(x)}$ croît strictement de $\abs{e_n(0)}=0$ à $\abs{e_n(1/2)}=2^{-n}\abs{E_n}$ (cf. \eqref{200307b}), lorsque $x$ croît de~$0$ à $1/2$, si $n$ est pair ;

$\bullet$ $\abs{e_n(x)}$ décroît strictement de $\abs{e_n(0)}$ à $\abs{e_n(1/2)}=0$, lorsque $x$ croît de $0$ à $1/2$, si $n$ est impair.

Pour déterminer la valeur de $e_n(0)$, considérons la série génératrice exponentielle de la suite~$E_n(0)$ :
\begin{align*}
\sum_{n\soe 0} E_n(0)\frac{z^n}{n!}&=\frac{2}{e^z+1} \expli{d'après \eqref{200307a}}\\
&=\frac{2}{e^z-1}-\frac{4}{e^{2z}-1}\\
&=2\sum_{n\soe 0} B_n\frac{z^{n-1}}{n!}-4\sum_{n\soe 0} B_n\frac{(2z)^{n-1}}{n!}\virg
\end{align*}
où les $B_n$ sont les nombres de Bernoulli, dont la série génératrice exponentielle est $z/(e^z-1)$. On a donc
\[
E_n(0)=\frac{B_{n+1}}{n+1}(2-2^{n+2}) \quad (n \in \Nat),
\]
et c'est aussi la valeur de $e_n(0)$, si $n \soe 1$. Ainsi, compte tenu des signes des nombres d'Euler et de Bernoulli, on a
\[
r_n=(-1)^{\floor{n/2}}\cdot\begin{dcases}
2^{-n}E_n & (n \text{ pair})\\
 (2^{n+2}-2)B_{n+1}/(n+1) &(n \text{ impair}).
\end{dcases}
\]
Les premières valeurs de la suite $(r_n)$ sont
\[
r_0=1 \quad ; \quad r_1=\demi\quad ; \quad r_2=\frac 14 \quad ; \quad r_3=\frac{1}{4} \quad ; \quad r_ 4=\frac{5}{16} \quad ; \quad r_ 5=\frac{1}{2}\cdotp
\]

\smallskip

Les \emph{splines d'Euler} sont définies par Schoenberg au moyen de la relation
\[
\Ecal_n(x) =e_n(x+\eps_n)/e_n(\eps_n) \quad (x \in \Real),
\]
où $\eps_n=0$ si $n$ est impair, et $\eps_n=1/2$ si $n$ est pair. On a donc $\abs{e_n(\eps_n)}=r_n$ et
\[
\norm{\Ecal_n}_{\infty}=\Ecal_n(0)=1.
\]

La fonction $\Ecal_n$ est une spline de degré $n$ sur $\Real$ ; ses principales propriétés ont été résumées par Schoenberg dans \cite{zbMATH03844309} ((VIII), p. 335) sous la forme suivante.
\begin{quote}
{\small La fonction $\Ecal_n(x)$ a toutes les propriétés de la fonction $\cos \pi x$ en ce qui concerne les symétries, les zéros, le signe et le sens de variation.}
\end{quote}

D'ailleurs, dans le même article, p. 336, Schoenberg explicite la série de Fourier de $\Ecal_n(x)$, et en déduit la formule asymptotique uniforme
\[
\Ecal_n(x) =\cos \pi x +O(3^{-n}).
\]

\smallskip

Dans les années précédant la parution de l'article de Kolmogorov \cite{zbMATH00050290}, les premières propriétés extrémales des fonctions $e_n$ furent découvertes. Pour $n \in \Nat^*$ et $\omega >0$, notons $X_n(\omega)$ l'ensemble des fonctions réelles de période $\omega>0$, $(n-1)$ fois dérivables, et à dérivée~$(n-1)$\up{e} lipschitzienne de rapport $1$.

Bernstein démontra en 1935 (cf. la version russe et corrigée \cite{zbMATH03089292}, p. 170-172, de la note \cite{zbMATH03017788}, et l'article de Favard \cite{zbMATH03025135} de 1936), que la borne supérieure de $\norm{f}_{\infty}$, lorsque $f \in X_n(\omega)$ est de moyenne nulle, était $(\omega/2)^nr_n/n!$, borne atteinte seulement pour les translatées de la fonction 
\begin{equation}\label{200317a}
f_n(x)=(\omega/2)^n\,\frac{e_n(2x/\omega)}{n!} \quad ( x \in \Real).
\end{equation}
 
  Si $\omega=2\pi$, ces valeurs extrémales sont parfois notées $K_n$ et appelées \emph{constantes de Favard}. En développant les fonctions $e_n$ en séries de Fourier, on peut montrer que
\[
K_n=\frac{4}{\pi}\cdot
\begin{dcases}
\sum_{k\soe 0}\,(-1)^k(2k+1)^{-n-1} & (n \text{ pair})\\
\sum_{k\soe 0}\,(2k+1)^{-n-1} & (n \text{ impair})
\end{dcases}
\]
Les deux suites $(K_{2m})_{m \in \Nat}$ et $(K_{2m+1})_{m \in \Nat}$ sont adjacentes ; leur limite commune est $4/\pi$.

Une autre propriété, découverte en 1936 par Favard (cf. la note \cite{zbMATH03025136} et l'article \cite{zbMATH03027646} de 1937 ; cf. également l'article d'Akhiezer et Krein \cite{zbMATH03026082} de 1937) est le fait que la fonction $f=f_n$, définie par \eqref{200317a}, réalise le maximum de la distance
\[
E_m(f)=\inf_{p \in \Tcal_{m-1}(\omega)} \norm{f-p}_{\infty},
\]
lorsque $f$ décrit $X_n(\omega)$, où $\Tcal_{m-1}(\omega)$ désigne l'ensemble des polynômes trigonométriques de période $\omega$ et de degré $<m$. Pour $m=1$, cette propriété découle de celle démontrée par Bernstein. On a $E_m(f_n)=(\omega/2)^nr_n/m^nn!$.

\subsection{Le théorème de Kolmogorov}\label{200416d}

Nous avons vu que $\|\Ecal_n\|_{\infty}=1$ pour tout $n$. Comme, pour presque tout $x$, on a 
\[
\Ecal_n^{(n)}(x)=e_n^{(n)}(x+\eps_n)/e_n(\eps_n)=n!e_0(x+\eps_n)/e_n(\eps_n),
\]
on a $\|\Ecal_n^{n)}\|_{\infty}=n!/r_n$. Notons $s_n=r_n/n!$ et
\[
q_n(x)=\Ecal_n(xs_n^{1/n}) \quad (x \in \Real),
\]
de sorte que $q_n \in \Lcal_n(1,1 \, ; \, \Real)$. La fonction $q$ du \S\ref{200117b} n'est autre que $q_2$.

Pour les dérivées intermédiaires de $q_n$, on a
\begin{equation*}
\|q_n^{(k)}\|_{\infty}=s_n^{k/n}\|\Ecal_{n}^{(k)}\|_{\infty}=\frac{s_n^{k/n}}{r_n}\|e_{n}^{(k)}\|_{\infty}
=\frac{n!s_n^{k/n}}{(n-k)!r_n}\|e_{n-k}\|_{\infty}=\frac{n!s_n^{k/n}r_{n-k}}{(n-k)!r_n}=s_{n-k}/s_n^{1-k/n},
\end{equation*}
d'après \eqref{200317b}.

Le fait remarquable, conjecturé par Chilov, et démontré par Kolmogorov, est que $\|q_n^{(k)}\|_{\infty}$ est le maximum de $\|f^{(k)}\|_{\infty}$, lorsque $f$ décrit $\Lcal_n(1,1; \, \Real)$, pour tout $k$ tel que~$0\ioe k\ioe n$. Par les considérations d'homogénéité du \S\ref{190322a}, on en déduit l'énoncé suivant.
\begin{theo}[Kolmogorov 1939, \cite{zbMATH00050290}]
Soit $n$ un nombre entier supérieur ou égal à $2$. Si $a$ et~$b$ sont deux nombres réels positifs, on a, pour toute fonction $f$ appartenant à $\Lcal_n(a,b; \, \Real)$,
\[
\|f^{(k)}\|_{\infty} \ioe \frac{s_{n-k}}{s_n^{1-k/n}} \,a^{1-k/n}\cdot b^{k/n}\quad (0 \ioe k \ioe n).
\]

Ces inégalités deviennent des égalités pour
\[
f(x)=aq_n\big(x(b/a)^{1/n}\big) \quad ( x \in \Real).
\]
\end{theo}

Notons que $s_{n-k}/s_n^{1-k/n}=K_{n-k}/K_n^{1-k/n}$. Ces quantités sont majorées, indépendamment de~$n$ et $k$, par $K_1=\pi/2$.

La démonstration de Kolmogorov repose sur le fait essentiel suivant : la fonction $q_n$ est une fonction de comparaison pour $\Lcal_n(1,1; \, \Real)$, \cad que
\[
\forall\, f \in \Lcal_n(1,1\, ; \, \Real), \; \forall \,t_0,t_1 \in \Real, \; f(t_0)=q_n(t_1) \implique \abs{f'(t_0)} \ioe \abs{q'_n(t_1)}.
\]

Bang (1941, cf. \cite{zbMATH03040318}) a donné une démonstration du théorème procédant en deux temps : on traite d'abord le cas des fonctions $f$ périodiques (et même presque périodiques), puis on approche une fonction quelconque de~$\Lcal_n(a,b; \, \Real)$ par une fonction périodique. La même démarche est suivie par Cavaretta dans \cite{zbMATH03449182}.

Par ailleurs, Bang donne une formulation du théorème en termes d'intégrales (plutôt que de dérivées) successives. Soit $f_0$ une fonction mesurable bornée sur $\Real$ et $(f_n)_{n \in \Nat}$ une suite obtenue par intégrations successives,
\[
f_n(x)=\int_0^x f_{n-1}(t) \, dt +d_n \quad (n \in \Nat^*, \; x \in \Real),
\]
où $(d_n)_{n \in \Nat}$ est une suite réelle arbitraire. Posons $M_n=\norm{f_n}_{\infty}$. Seule $M_0$ est, a priori, supposée finie, et non nulle. Le théorème de Kolmogorov, vu par Bang (cf. \cite{zbMATH03040318}, p. 7), est alors l'assertion selon laquelle la suite $\rho_n=(M_n/M_0s_n)^{1/n}$ (pour $n \soe 1$) est croissante au sens large. Lorsque l'on prend $f_n=e_n/n!$  pour tout $n$, la suite $\rho_n$ correspondante est la constante $1$.

\subsection{Le cas $I=[0,T]$}\label{200202a}

Un résultat analogue à \eqref{200330a} n'est connu explicitement que pour une seule autre valeur de $n$, à savoir $n=3$. Voici ce résultat, dû à M. Sato (1982, cf. \cite{zbMATH03778028}) :
\begin{align*}
\sup_{f \in \Lcal_3(a,b \, ;\, [0,T])} \norm{f'}_{\infty} &=
\begin{dcases}
4 \cdot a/(\alpha T) + \frac 16 \cdot b(\alpha T)^2 & ( T \ioe T_0)\\
C_{3,1}\cdot a^{2/3}b^{1/3} & ( T \soe  T_0)
\end{dcases}
\\
\sup_{f \in \Lcal_3(a,b \, ;\, [0,T])} \norm{f''}_{\infty} &=
\begin{dcases}
4\cdot a/(\alpha T)^2 + \frac 23 \cdot b(\alpha T )& ( T \ioe  T_0)\\
C_{3,2} \cdot a^{1/3}\,b^{2/3} & ( T \soe  T_0),
\end{dcases}
\end{align*}
où 
\begin{equation}\label{200412a}
T_0=T_0(a,b)=(81a/b)^{1/3} \quad ; \quad
C_{3,1} =3^{5/3}/2 \quad ; \quad
C_{3,2} =2 \cdot 3^{1/3},
\end{equation}
et où $\alpha$ est la solution dans  $[1/3,1/2[$ de l'équation biquadratique
\[ 
12-24 \alpha=\frac{T^3b}a \alpha^2(1-\alpha)^2.
\]

\smallskip

Ce résultat a notamment l'intérêt de mettre en évidence que, pour $T\ioe T_0$, les bornes supérieures cherchées ne sont pas de la forme du second membre de l'inégalité \eqref{200402b}, avec des constantes~$A_{3,k}$ et $B_{3,k}$ indépendantes de $T$.

Cela étant, dans le cas général, où $n$ est un nombre entier quelconque, plusieurs auteurs ont explicité des constantes~$A_{n,k}$ et $B_{n,k}$ admissibles dans la proposition~\ref{200402a}. Une première observation est que l'inégalité \eqref{200402b} s'applique en particulier au cas où $f=P$ est une fonction polynomiale de degré au plus $n-1$, vérifiant $\abs{f} \ioe 1$ sur $[0,2]$ (ou $[-1,1]$). On peut alors prendre $a=1$, $T=2$, et $b$ arbitrairement petit. Ainsi,
\[
A_{n,k} \soe 2^k\sup_{\deg P \ioe n-1}\max_{\abs{t}\ioe 1}\, \lvert P^{(k)}(t)\rvert.
\]

Or cette borne supérieure est connue, grâce à un théorème des frères Markov (cf. le survol de A. Shadrin \cite{MR2117091}) : elle est atteinte lorsque $P$ est le $(n-1)$\up{e} polynôme de Tchebychev, défini par l'identité $T_{n-1}(x)=\cos (n-1) \arccos x$, et pour $t=1$. On a donc
\[
A_{n,k} \soe 2^k \lvert T_{n-1}^{(k)}(1)\rvert.
\]

Notons que
\[
\lvert T_{n-1}^{(k)}(1)\rvert=T_{n-1}^{(k)}(1)=(n-1)\frac{ 2^k k!}{(2k)!} \frac{(n+k-2)!}{(n-k-1)!}\cdotp
\]

En approchant une fonction $f\in \Lcal_n(a,b \, ; [0,T])$ par son polynôme de meilleure approximation de degré~$\ioe n-1$, Gorny avait obtenu en 1939 une valeur $A_{n,k}$ double de cette borne inférieure (cf. \cite{zbMATH03035809}, (3), p. 321). L'année suivante, Cartan fut le premier a obtenir l'inégalité \eqref{200402b} avec la valeur optimale, $A_{n,k} = 2^k T_{n-1}^{(k)}(1)$.

Une démarche possible est alors de choisir cette valeur pour $A_{n,k}$, et de chercher une constante admissible~$B_{n,k}$ la plus petite possible. En prenant $T=2$, $f(t)=T_n(t-1)$ dans \eqref{200402b}, on obtient
\[
T_{n}^{(k)}(1) \ioe T_{n-1}^{(k)}(1)+B_{n,k}T_{n}^{(n)}(1)2^{n-k},
\]
d'où
\begin{equation}\label{200411a}
B_{n,k} \soe A_{n,k} \cdot \frac{k(2n-1)}{(n-k)(n-1)\, n! 2^{2n-1}}\cdotp
\end{equation}

En 1940, Cartan obtenait la valeur admissible
\[
B_{n,k}=A_{n,k}/n!
\]
(cf. \cite{zbMATH03099983}, (7), p. 11, et (12), p. 13). En 1990, Kallioniemi a obtenu la valeur
\begin{equation}\label{200412b}
B_{n,k} = A_{n,k} \cdot \frac{k\big(2(n-1)^2+k-1\big)}{(n-k)(n-1)(n+k-2)\, n! 2^{2n-2}}
\end{equation}
(cf. \cite{zbMATH04214632}, Theorem 4.3, p. 84), dont le quotient par le second membre de \eqref{200411a} est compris entre $1$ et $2$.

\smallskip

La recherche actuelle sur ce problème est structurée par une conjecture de 1976, due à Karlin (cf. \cite{zbMATH03504851}, p. 423), et stipulant qu'une certaine fonction spline $Z_T$, dite \emph{spline de Zolotareff}, et précisément décrite par Karlin\footnote{Karlin donne cette description pour l'ensemble $\Lcal_n(\rho,1  \, ; [0,1])$, où $\rho>0$ ; le passage à notre $\Lcal_n(1,1  \, ; \,[0,T])$ s'effectuerait suivant les principes d'homogénéité du \S\ref{190322a}.} (\cite{zbMATH03504851}, Theorem 5.1, p. 376, et \S 5, p. 403-411), vérifie
\[
\sup_{f \in \Lcal_n(1,1 \, ;\, [0,T])}\|f^{(k)}\|_{\infty} = \|Z_T^{(k)}\|_{\infty} =\lvert Z_T^{(k)}(0)\rvert \quad (0 <k<n).
\]

\smallskip

En particulier, lorsque $T \ioe 4(n!/2)^{1/n}$, la fonction spline $Z_T$ est un polynôme, obtenu par changement de variable entre $[-1,1]$ et $[0,T]$, à partir d'un \emph{polynôme de Zolotareff}. Ces polynômes ont été introduits en 1868 par Zolotareff dans sa thèse d'habilitation (cf. \cite{zbMATH03008070}, p.~130-166), dont voici une traduction française des premières lignes.
\begin{quote}
{\small
Le problème, dont je donne la solution dans cet article, est le suivant :

\smallskip

{\og De la fonction entière (polynomiale)
\[
F(x)=x^n-\sigma x^{n-1} + \cdots,
\]
où $\sigma$ est un coefficient fixé, trouver les coefficients restants de sorte que cette fonction garde la plus petite valeur possible lorsque $x$ reste entre les bornes $-1$ et $1$, et trouver le maximum de cette fonction.\fg}

\smallskip

Pour $\sigma=0$, cette question a été résolue par P. L. Tchebychev dans le mémoire \emph{Théorie des mécanismes connus sous le nom de parallélogrammes}, 1853.

Dans ce cas, $F(x)$ s'exprime très simplement à l'aide des fonctions trigonométriques ; elle est  précisément égale à $2^{1-n}\cos n \arccos x$.

Lorsque $\sigma$ n'est pas nul, nous démontrons que $F(x)$ s'exprime très simplement à l'aide des fonctions de Jacobi (...)
}
\end{quote}

\smallskip

Pour le cas polynomial ($T \ioe 4(n!/2)^{1/n}$), la conjecture de Karlin a été démontrée en 2014 par Shadrin (cf. \cite{zbMATH06224681}).

\subsection{Le cas $I=[0,\infty[$}\label{200416e}

Posons
\[
C_{n,k}=\sup_{f \in \Lcal_n(1,1; \, [0,\infty[)}\|f^{(k)}\|_{\infty} \quad (n \soe 2\, ; \, 0<k<n).
\]

L'égalité $C_{2,1}=2$, due à Landau, est l'objet de la proposition \ref{190627g} ci-dessus. Les valeurs $C_{3,1}$ et $C_{3,2}$ de \eqref{200412a}, ont été obtenues par Matorin en 1955 (cf. \cite{matorin}).

Schoenberg et Cavaretta ont montré en 1970 (cf. \cite{zbMATH03364147}) l'existence d'une suite $(q_n^+)_{n \soe 2}$ de fonctions vérifiant

\smallskip

$\bullet$ $q_n^+$ est une spline de degré $n$ sur $[0,\infty[$ ;

$\bullet$ $q_n^+ \in \Lcal_n(1,1; \, [0,\infty[)$ pour tout $n \soe 2$ ;

$\bullet$ $\norm{f^{(k)}}_{\infty} \ioe \|q_n^{+(k)}\|_{\infty}$ pour $n \soe 2$, $f \in \Lcal_n(1,1; \, [0,\infty[)$, et $0<k<n$.

\smallskip

Les fonctions $q_n^+$ sont explicitement connues pour $n=2$ et $3$ (cf. \cite{zbMATH03411359}, pp. 147-156). Pour~$n \soe 4$, les fonctions $q_n^+$ sont définies par un passage à la limite, donnant des estimations des bornes~$C_{n,k}$ qui sont {\og far from explicit\fg} (cf.~\cite{zbMATH03364147}, p. 304). 

\smallskip

Concernant le comportement asymptotique de ces constantes, la question principale est la suivante.

\begin{ques}\label{200416g}
Soit $\lambda$ tel que $0 < \lambda <1$. La quantité $C_{n,k}^{1/n}$ tend-elle vers une limite quand $n$ tend vers l'infini et $k/n$ tend vers $\lambda$? Si la réponse est oui, que vaut cette limite?
\end{ques}

Pour $0 < \lambda <1$, posons
\[
w^+(\lambda)=\limsup \frac{\ln C_{n,k}}{n} \quad ; \quad w^-(\lambda)=\liminf \frac{\ln C_{n,k}}{n} \quad (n \vers \infty, \; k/n \vers \lambda).
\]

En 1955, en utilisant la méthode de comparaison de Kolmogorov, Matorin obtint la majoration
\[
C_{n,k} \ioe T_n^{(k)}(1)/T_n^{(n)}(1)^{k/n}
\]
(cf. \cite{matorin}, (3), p. 13). La formule de Stirling sous la forme
\[
(k!)^{1/n}\sim \lambda^{\lambda}e^{-\lambda}n^{\lambda} \quad (n \vers \infty, \; k/n \vers \lambda)
\]
permet d'en déduire la majoration
\begin{equation}\label{200412c}
w^+(\lambda) \ioe (1+\lambda)\ln (1+\lambda)+(\lambda -1)\ln (1- \lambda)-\lambda \ln (4 \lambda) \quad (0 < \lambda <1).
\end{equation}

La même année, par une méthode entièrement différente, fondée sur des considérations d'analyse complexe et d'analyse fonctionnelle, Malliavin obtenait dans sa thèse la majoration
\begin{equation}\label{200415a}
C_{n,k} \ioe 2^{10}\frac{e\ln n}{\pi} \, e^{n\mu(k/n)} \quad \big(\,\mu(\lambda)=-\int_0^{\lambda} \ln \tan (\pi t/2) \, dt\, \big)
\end{equation}
(cf. \cite{zbMATH03113668}, 11.3.4, p. 238, et 14.6, p. 254), d'où découle immédiatement l'inégalité $w^+(\lambda)\ioe \mu(\lambda)$. On peut vérifier que cette majoration est toujours meilleure que \eqref{200412c}.

\smallskip

En ce qui concerne $w^-(\lambda)$, la minoration obtenue en 1967 par Stechkin,
\[
C_{n,k} \soe \kappa p^{-1/2}(n/p)^p, 
\]
où $p=\min(k,n-k)$, et où $\kappa$ est une constante positive absolue (cf. \cite{stechkin}, (23), p. 445), fournit l'inégalité
\[
w^-(\lambda) \soe \min(\lambda,1-\lambda) \ln 1/\min(\lambda,1-\lambda) \quad ( 0 < \lambda <1).
\]

\subsection{Autres normes}\label{200416f}

Pour conclure, voici quelques indications sur un problème généralisant ceux considérés dans ce survol. On se donne, comme précédemment, un intervalle $I$ de $\Real$, deux nombres entiers $k$ et $n$, vérifiant l'encadrement~$0<k<n$, deux nombres réels positifs $a$ et $b$, mais aussi trois éléments~$p,\, q , \, r$ de~$[1,\infty]$, et on se propose de déterminer la borne supérieure 
\[
\sigma_{p,q,r}(a,b\, ; \, k,n \, ; \, I)=\sup\|f^{(k)}\|_q ,
\]
lorsque $f$ parcourt l'ensemble des fonctions $f : I \vers \Real$, qui sont $(n-1)$ fois dérivables, dont la dérivée~$(n-1)$\up{e} est localement absolument continue dans $I$, et telles que
\[
\|f\|_p \ioe a \quad ; \quad \|f^{(n)}\|_r \ioe b.
\]

L'essentiel de ce texte était dévolu au cas $p=q=r=\infty$, mais nous avons vu, par exemple, l'égalité
\[
\sigma_{\infty,1,\infty}(1,1\, ; \, 1,2 \, ; \, [0,2N\sqrt{2}\, + 4])=2N+4 \quad ( N \in \Nat)
\]
(proposition \ref{190615b} ci-dessus). Notons que, dans l'article \cite{zbMATH01891483}, dont ce résultat est issu, Bojanov et Naidenov obtiennent une égalité similaire pour $n=2$ ou $3$, $p=r=\infty$, et $q$ quelconque.

Le livre \cite{zbMATH00205893} de Kwong et Zettl contient une introduction au problème général.

\smallskip

Si $I=\Real$, Stein a montré en 1957 que l'inégalité
\[
\|f^{(k)}\|_{p} \ioe \frac{s_{n-k}}{s_n^{1-k/n}} \cdot\|f\|_{p}^{1-k/n}\cdot \|f^{(n)}\|_{p}^{k/n}\quad (0 \ioe k \ioe n),
\]
due à Kolmogorov pour $p=\infty$, est valable pour toute valeur de $p \in [1,\infty]$ (cf. \cite{zbMATH03129755}, p. 586-588). Outre $p=\infty$, cette majoration est aussi optimale pour $p=1$ (Ditzian, 1975, cf. \cite{zbMATH03484694}, p. 149). Pour $p=2$, la majoration optimale est simplement
\[
\|f^{(k)}\|_{2} \ioe \|f\|_{2}^{1-k/n}\cdot \|f^{(n)}\|_{2}^{k/n}\quad (0 \ioe k \ioe n),
\]
et $p=1,2,\infty$ sont les seules valeurs pour lesquelles les bornes $\sigma_{p,p,p}(a,b\, ; \, k,n \, ; \, \Real)$ soient connues.

\smallskip

Si $I=[0,\infty[$, le premier résultat est l'inégalité démontrée en 1932 par Hardy et Littlewood,
\[
\|f'\|_{2} \ioe \sqrt{2} \;\|f\|_{2}^{1/2}\cdot \|f''\|_{2}^{1/2},
\]
où la constante $\sqrt{2}$ est optimale (cf. \cite{zbMATH03007670}, Theorem 6, p. 249). Dans \cite{zbMATH03113668} (14.1, p. 251), Malliavin énonce l'inégalité générale
\[
\|f^{(k)}\|_{2} \ioe e^{n\mu(k/n)} \;\|f\|_{2}^{1-k/n}\cdot \|f^{(n)}\|_{2}^{k/n} \quad (0 <k <n)
\]
(où $\mu$ est définie dans \eqref{200415a}), dans le cas où $f$ est indéfiniment dérivable sur $]0,\infty[$.

\smallskip

Enfin, si $I=[0,T]$, les résultats actuellement les plus précis sont ceux des articles de Bojanov et Naidenov déjà cités.

\begin{center}
{\sc Remerciements}
\end{center}

{\footnotesize  Je remercie Sonia Fourati de m'avoir invité et encouragé à écrire cet article, et l'arbitre anonyme d'avoir suggéré l'écriture du \S\ref{200311a}.}


\medskip

\footnotesize

\noindent BALAZARD, Michel\\
Aix Marseille Univ, CNRS, Centrale Marseille, I2M, Marseille, France\\
Adresse \'electronique : \texttt{balazard@math.cnrs.fr}

\end{document}